\documentclass[11pt,reqno]{amsart}
\usepackage{geometry}
\usepackage{graphicx}
\usepackage{amssymb}
\usepackage{caption}
\usepackage{subcaption}
\usepackage{enumitem}
\usepackage{mathtools}
\usepackage[all]{xy}
\usepackage{xcolor}
\usepackage{soul}
\usepackage[hidelinks,colorlinks=true,linkcolor=blue,urlcolor=blue,citecolor=blue]{hyperref}

\usepackage{amsmath}

\usepackage{MnSymbol}

 \usepackage{dutchcal}

\usepackage[linesnumbered,ruled,vlined]{algorithm2e}

\newcommand{\noln}{\renewcommand{\nl}{\let\nl\oldnl}}

\SetCommentSty{mycommfont}
\SetKwInput{KwInput}{Input}               
\SetKwInput{KwOutput}{Output}

\graphicspath{{figures/}}

\numberwithin{equation}{section}

\newtheorem{theorem}{Theorem}[section]

\newtheorem{remark}{Remark}[section]
\newtheorem*{remark*}{Remark}

\newtheorem{definition}{Definition}[section]

\usepackage[numbers,sort&compress]{natbib}
\bibpunct{[}{]}{;}{n}{,}{,}

\DeclarePairedDelimiterX{\inner}[2]{\langle}{\rangle}{#1, #2}

\DeclareMathOperator*{\ext}{ext}

\allowdisplaybreaks

\title[Accelerated Optimization via Discrete Constrained Variational Integrators]{Accelerated Optimization on Riemannian Manifolds via 
	\\ Discrete Constrained Variational Integrators}
\author{Valentin Duruisseaux}
\author{Melvin Leok}
\thanks{ Corresponding author: Melvin Leok, \texttt{mleok@ucsd.edu}}
\address{Department of Mathematics, University of California, San Diego, La Jolla, CA 92093-0112, USA.}
\email{vduruiss@ucsd.edu, mleok@ucsd.edu}

\begin{document}

\maketitle

\begin{abstract}
	A variational formulation for accelerated optimization on normed vector spaces was recently introduced in~\citet{WiWiJo16}, and later generalized to the Riemannian manifold setting in~\citet{Duruisseaux2021Riemannian}. This variational framework was exploited on normed vector spaces in~\citet{duruisseaux2020adaptive} using time-adaptive geometric integrators to design efficient explicit algorithms for symplectic accelerated optimization, and it was observed that geometric discretizations which respect the time-rescaling invariance and symplecticity of the Lagrangian and Hamiltonian flows were substantially less prone to stability issues, and were therefore more robust, reliable, and computationally efficient. As such, it is natural to develop time-adaptive Hamiltonian variational integrators for accelerated optimization on Riemannian manifolds. In this paper, we consider the case of Riemannian manifolds embedded in a Euclidean space that can be characterized as the level set of a submersion. We will explore how holonomic constraints can be incorporated in discrete variational integrators to constrain the numerical discretization of the Riemannian Hamiltonian system to the Riemannian manifold, and we will test the performance of the resulting algorithms by solving eigenvalue and Procrustes problems formulated as optimization problems on the unit sphere and Stiefel manifold.
\end{abstract}

\section{Introduction}

Many data analysis algorithms are designed around the minimization of a loss function or the maximization of a likelihood function. Due to the ever-growing scale of the data sets and size of the problems, there has been a lot of focus on first-order optimization algorithms because of their low cost per iteration. In 1983, Nesterov's accelerated gradient method~\cite{Nes83} was shown to converge in $\mathcal{O}(1/k^2)$ to the minimum of the convex objective function $f$, improving on the $\mathcal{O}(1/k)$ convergence rate exhibited by standard gradient descent methods. This $\mathcal{O}(1/k^2)$ convergence rate was shown in~\cite{Nes04} to be optimal among first-order methods using only information about $\nabla f$ at consecutive iterates. This phenomenon in which an algorithm displays this improved rate of convergence is referred to as acceleration, and other accelerated algorithms have been derived since Nesterov's algorithm, which was shown in~\cite{SuBoCa16} to limit to a second-order ordinary differential equation (ODE), as the timestep goes to 0, and that $f(x(t))$ converges to its optimal value at a rate of $\mathcal{O}(1/t^2)$ along any trajectory $x(t)$ of this ODE. It was then shown in~\cite{WiWiJo16} that in continuous time, an arbitrary convergence rate $\mathcal{O}(1/t^p)$ can be achieved in normed spaces, by considering flow maps generated by a family of time-dependent Bregman Lagrangian and Hamiltonian systems which is closed under time-rescaling. This variational framework and the time-rescaling property of this family were then exploited in~\cite{duruisseaux2020adaptive} by using time-adaptive geometric integrators to design efficient explicit algorithms for symplectic accelerated optimization. It was observed that a careful use of adaptivity and symplecticity could result in a significant gain in computational efficiency. More generally, when applied to Hamiltonian systems, symplectic integrators yield discrete approximations of the flow that preserve the symplectic 2-form~\cite{HaLuWa2006}. The preservation of symplecticity results in the preservation of many qualitative aspects of the underlying dynamical system. In particular, when applied to conservative Hamiltonian systems, symplectic integrators exhibit excellent long-time near-energy preservation~\cite{Benettin1994,Re1999}. Variational integrators provide a systematic method for constructing symplectic integrators of arbitrarily high-order based on the numerical discretization of Hamilton's principle~\cite{MaWe2001, HaLe2012}, or equivalently, by the approximation of Jacobi's solution of the Hamilton--Jacobi equation, which is a generating function for the exact symplectic flow map.

In the past few years, there has been some effort to derive accelerated optimization algorithms in the Riemannian manifold setting~\cite{Duruisseaux2021Riemannian,alimisis2020,Alimisis2021,Alimisis2020-1,Sra2016, Sra2018, Sra2020, Liu2017}. In~\cite{Duruisseaux2021Riemannian}, it was shown that in continuous time, the convergence rate of $f(x(t))$ to its optimal value can be accelerated to an arbitrary convergence rate $\mathcal{O}(1/t^p)$ on Riemannian manifolds. This was achieved by considering a family of time-dependent Bregman Lagrangian and Hamiltonian systems on Riemannian manifolds which is closed under time-rescaling, thereby generalizing the variational framework for accelerated optimization of~\cite{WiWiJo16} to Riemannian manifolds. The time-adaptivity based approach relying on a Poincar\'e transformation from~\cite{duruisseaux2020adaptive} was also extended to the Riemannian manifold setting in~\cite{Duruisseaux2021Riemannian}. Now, the Whitney Embedding Theorems~\cite{Whitney1944_2,Whitney1944_1} state that any smooth manifold of dimension~$n \geq 2$ can be embedded in $\mathbb{R}^{2n}$ and immersed in $\mathbb{R}^{2n-1}$, and is thus diffeomorphic to a submanifold of~$\mathbb{R}^{2n}$. Furthermore, the Nash Embedding Theorems~\cite{Nash1956} state that any Riemannian manifold can be globally isometrically embedded into some Euclidean space. As a consequence of these embedding theorems, the study of Riemannian manifolds can in principle be reduced to the study of submanifolds of Euclidean spaces. Altogether, this motivates the introduction of time-adaptive variational integrators on Riemannian manifolds which exploit the structure of the embedding Euclidean space, and in this paper we will study how holonomic constraints can be incorporated into different types of variational integrators to constrain the numerical solutions of the Riemannian dynamical system to the Riemannian manifold. Incorporating holonomic constraints in geometric integrators has been studied extensively in the past (see~\cite{MaWe2001, HaLuWa2006, MaRa1999, Holm2009} for instance), and some work has been done from the variational perspective for the Type I Lagrangian formulation in~\cite{MaWe2001} via augmented Lagrangians. 

\subsection*{Outline of the paper}

We will first show in Section~\ref{section: Constrained Variational Lagrangian} the equivalence between constrained variational principles and constrained Euler--Lagrange equations, both in continuous and discrete time, before deriving analogous results for both the Type II and Type III Hamiltonian formulations of mechanics in Section~\ref{Sec: Constrained H Mechanics}. In Section~\ref{sec: Error Analysis}, we will exploit error analysis theorems for unconstrained mechanics from~\cite{MaWe2001,ScLe2017} to obtain variational error analysis results for the maps defined implicitly by the discrete constrained Euler--Lagrange and Hamilton's equations. Finally, in Section~\ref{sec: Optimization Section}, we will exploit these constrained variational integrators and the variational formulation of accelerated optimization on Riemannian manifolds from~\cite{Duruisseaux2021Riemannian} to solve numerically generalized eigenvalue problems and Procrustes problems on the unit sphere and Stiefel manifold.

\section{Constrained Variational Lagrangian Mechanics}  \label{section: Constrained Variational Lagrangian}

Traditionally, variational integrators have been designed based on the Type I generating function known as the \textbf{discrete Lagrangian}, $L_d:\mathcal{Q} \times \mathcal{Q} \rightarrow \mathbb{R}$. The exact discrete Lagrangian is the exact generating function for the time-$h$ flow map of Hamilton's equations and it can be represented in boundary-value form by
\begin{equation}
	L_d^E(q_0,q_h)=\int_0^h L(q(t),\dot q(t)) dt  \label{exact_Ld},
\end{equation}
where $q(0)=q_0,$ $q(h)=q_h,$ and $q$ satisfies the Euler--Lagrange equations over the time interval $[0,h]$. This is closely related to Jacobi's solution of the Hamilton--Jacobi equation. A variational integrator is defined by constructing an approximation $L_d:\mathcal{Q} \times \mathcal{Q} \rightarrow \mathbb{R}$ to the exact discrete Lagrangian $L_d^E$, and then applying the \textbf{implicit discrete Euler--Lagrange equations},
\begin{equation}
	p_0=-D_1 L_d(q_0, q_1),\qquad p_1=D_2 L_d(q_0, q_1),  \label{IDEL}
\end{equation}
which implicitly define a numerical integrator, referred to as the \textbf{discrete Hamiltonian map} $\tilde{F}_{L_d}:(q_0,p_0)\mapsto(q_1,p_1)$, where $D_i$ denotes a partial derivative with respect to the $i$-th argument. These equations define the \textbf{discrete Legendre transforms}, \(\mathbb{F}^{\pm}L_{d}: \mathcal{Q}\times \mathcal{Q} \rightarrow T^{*}\mathcal{Q}\):
\begin{align}
\mathbb{F}^{+}L_{d}&:(q_{0},q_{1}) \mapsto (q_{1},p_{1}) = (q_{1},D_{2}L_{d}(q_{0},q_{1})), \\
\mathbb{F}^{-}L_{d}&:(q_{0},q_{1}) \mapsto (q_{0},p_{0}) = (q_{0},-D_{1}L_{d}(q_{0},q_{1})),
\end{align}
and the discrete Hamiltonian map can be expressed as $\tilde{F}_{L_d}\equiv (\mathbb{F}^{+}L_{d})\circ (\mathbb{F}^{-}L_{d})^{-1}$. Such numerical methods are called variational integrators as they can be derived from a \textbf{discrete Hamilton's principle}, which involves extremizing a discrete action sum $S_d\left(\{q_k\}_{k=0}^N \right)\equiv \sum_{k=0}^{N-1} L_d(q_k,q_{k+1})$, subject to fixed boundary conditions on $q_0$, $q_N$. 

Now, suppose we are given a configuration manifold $\mathcal{M}$, and a holonomic constraint function $\mathcal{C} : \mathcal{M} \rightarrow \mathbb{R}^d$. Assuming that $0\in \mathbb{R}^d$ is a regular point of $\mathcal{C}$, we can constrain the dynamics to the constraint submanifold $\mathcal{Q} = \mathcal{C}^{-1}(0),$ which is truly a submanifold of $\mathcal{M}$ (see~\cite{MaWe2001,AbMaRa1988}). We will now consider variational Lagrangian mechanics with holonomic constraints $\mathcal{C}(q)$ using Lagrange multipliers $\lambda : [0,T] \rightarrow \Lambda$.

\subsection{Continuous Constrained Variational Lagrangian Mechanics}

We begin by presenting an equivalence between the continuous constrained variational principle and the continuous constrained Euler--Lagrange equations:

\begin{theorem} \label{theorem: Lagrangian Action and Equations}
	Consider the \textbf{constrained action functional} $\mathfrak{S} : C^2([0,T], \mathcal{Q} \times \Lambda) \rightarrow \mathbb{R}$ given by
	\begin{equation}\label{Lagrangian Constrained Action Functional}
		\mathfrak{S} (q(\cdot), \lambda(\cdot)) =  \int_{0}^{T}{ \left[ L(q(t),\dot{q}(t) ) - \langle \lambda(t) , \mathcal{C}(q(t)) \rangle \right] dt}.
	\end{equation}
	The condition that $\mathfrak{S} (q(\cdot), \lambda(\cdot)) $ is stationary with respect to the boundary conditions $\delta q(0) = 0$ and $\delta q(T) =0$ is equivalent to $\left(q(\cdot),\lambda(\cdot) \right)$ satisfying the \textbf{constrained Euler--Lagrange equations}
	\begin{equation} \label{ELConstEqs}
			\frac{\partial L}{\partial q} - \frac{d}{dt} \frac{\partial L}{\partial \dot{q}}  =  \langle \lambda , \nabla \mathcal{C} (q) \rangle, \qquad \mathcal{C}(q) =0 . 
	\end{equation}
	\proof{See Appendix~\ref{Appendix: Type I Action}.  \qed}
\end{theorem}

\begin{remark}
	These constrained Euler--Lagrange equations can be thought of as the Euler--Lagrange equations coming from the augmented Lagrangian $\bar{L} \left( q,\lambda, \dot{q}, \dot{\lambda} \right)  = L(q,\dot{q} )- \langle \lambda , \mathcal{C} (q) \rangle$.
\end{remark}

Consider the function $\mathcal{S}(q_0,q_T)$ given by the extremal value of the constrained action functional~$\mathfrak{S}$ over the family of curves $(q(\cdot) , \lambda(\cdot)) $ satisfying the boundary conditions $q(0)=q_0$ and $q(T)=q_T$:
\begin{equation}
	\mathcal{S}(q_0,q_T) = \ext_{ \substack{(q,\lambda)\in C^2([0,T],\mathcal{Q} \times \Lambda)   \\ q(0)=q_0, \quad  q(T)=q_T}}{\mathfrak{S} (q(\cdot) , \lambda(\cdot)) }.
\end{equation}
The following theorem shows that $\mathcal{S}(q_0,q_T)$ is a generating function for the flow of the continuous constrained Euler--Lagrange equations:

\begin{theorem} \label{theorem: Flow Map and Equations I}
		The exact time-$T$ flow map of Hamilton's equations $(q_0,p_0) \mapsto (q_T,p_T)$ is implicitly given by the following relations:
	\begin{equation}
		D_1 \mathcal{S} (q_0,q_T) = -  \frac{\partial L}{ \partial \dot{q}}(q_0,\dot{q}(0)), \qquad  D_2 \mathcal{S} (q_0,q_T) =  \frac{\partial L}{ \partial \dot{q}}(q_T,\dot{q}(T))  .
	\end{equation}
	In particular, $ \mathcal{S} (q_0,q_T)$ is a Type I generating function that generates the exact flow of the constrained Euler--Lagrange equations \eqref{ELConstEqs}. 
	\proof{See Appendix~\ref{Appendix: Type I Flow}. \qed
	}
\end{theorem}

\subsection{Discrete Constrained Variational Lagrangian Mechanics} \label{sec: Discrete Lagrangian}

We now introduce a discrete variational formulation of Lagrangian  mechanics which includes holonomic constraints. Suppose we are given a partition $0 = t_0 < t_1 < \ldots < t_N = T$ of the interval $[0,T]$, and a discrete curve in $\mathcal{Q} \times \Lambda$ denoted by $\{  (q_k, \lambda _k) \}_{k=0}^{N}$ such that $q_k \approx q(t_k)$ and $\lambda_k \approx \lambda(t_k)$.
We will formulate discrete constrained variational Lagrangian mechanics in terms of the following discrete analogues of the constrained action functional $\mathfrak{S}$ given by equation \eqref{Lagrangian Constrained Action Functional}:
\begin{align} \label{Discrete Constrained Action Functional I}
		\mathfrak{S}_d^+ \left(\{  (q_k, \lambda _k) \}_{k=0}^{N} \right) & =  \sum_{k=0}^{N-1}{\left[ L_d(q_k,q_{k+1} )- \langle \lambda_{k+1} , \mathcal{C} (q_{k+1}) \rangle \right]}, \\
	\mathfrak{S}_d^- \left(\{  (q_k, \lambda _k) \}_{k=0}^{N} \right) & =  \sum_{k=0}^{N-1}{\left[ L_d(q_k,q_{k+1} )- \langle \lambda_k , \mathcal{C} (q_k) \rangle \right]},  
\end{align}
where
\begin{align}\label{Discrete Lagrangian}
	L_d(q_k,q_{k+1} )  & \approx  \ext_{ \substack{(q,\lambda)\in C^2([t_k,t_{k+1}],\mathcal{Q} \times \Lambda)   \\ q(t_{k})=q_{k}, \quad  q(t_{k+1})=q_{k+1}}}   \text{  } \int_{t_k}^{t_{k+1}}{L(q(t),\dot{q}(t)) dt}.
\end{align}

We can now derive discrete analogues to Theorem~\ref{theorem: Lagrangian Action and Equations} relating discrete Type I variational principles to discrete Euler--Lagrange equations:

\begin{theorem} \label{theorem: Discrete Variational Principle I } The Type I discrete Hamilton's variational principles 
	\begin{equation}
		\delta \mathfrak{S}_d^{\pm} \left(\{  (q_k, \lambda _k) \}_{k=0}^{N} \right)  = 0,
	\end{equation}
	are equivalent to the \textbf{discrete constrained Euler--Lagrange equations}
	\begin{equation}
		D_1 L_d(q_k,q_{k+1})  + D_2 L_d(q_{k-1},q_{k}) = \langle \lambda_k ,\nabla \mathcal{C}(q_k) \rangle, \quad  \mathcal{C}(q_{k})  =0,
	\end{equation}
	where $L_d(q_{k},q_{k+1}) $ is defined via equation \eqref{Discrete Lagrangian}. 
	\proof{See Appendix~\ref{Appendix: Type I Discrete Variational}. \qed
	}
\end{theorem}   

\begin{remark} \label{remark: augmented Lagrangian}
	These discrete constrained Euler--Lagrange equations can be thought of as the discrete Euler--Lagrange equations coming from the augmented discrete Lagrangians
\begin{align} 
		\bar{L}_d^+ \left( q_k,\lambda_k, q_{k+1}, \lambda_{k+1}\right) & = L_d(q_k,q_{k+1} )- \langle \lambda_{k+1} , \mathcal{C} (q_{k+1}) \rangle, \\
\bar{L}_d^- \left( q_k,\lambda_k, q_{k+1}, \lambda_{k+1}\right) & = L_d(q_k,q_{k+1} )- \langle \lambda_k , \mathcal{C} (q_k) \rangle. 
\end{align} \\
\end{remark}

\section{Constrained Variational Hamiltonian Mechanics} \label{Sec: Constrained H Mechanics}

The boundary-value formulation of the exact Type II generating function of the time-$h$ flow of Hamilton's equations is given by the exact discrete right Hamiltonian,
\begin{equation}
	H_d^{+,E}(q_0,p_h) =  p_h q_h - \int_0^h \left[ p(t) \dot{q}(t)-H(q(t), p(t)) \right] dt, \label{exact_Hd}
\end{equation}
where $(q,p)$ satisfies Hamilton's equations with boundary conditions $q(0)=q_0$ and $p(h)=p_h$. A Type II Hamiltonian variational integrator is constructed by using an approximate discrete right Hamiltonian $H_d^+$, and applying the \textbf{discrete right Hamilton's equations},
\begin{equation}\label{Discrete Right Eq}
	p_0=D_1H_d^+(q_0,p_1), \qquad q_1=D_2H_d^+(q_0,p_1),
\end{equation}
which implicitly define the integrator, $\tilde{F}_{H_d^+}:(q_0,p_0) \mapsto (q_1,p_1)$. 

Similarly, the boundary-value formulation of the exact Type III generating function of the time-$h$ flow of Hamilton's equations is given by the exact discrete left Hamiltonian,
\begin{align}
	H_d^{-,E}(q_h,p_0) =  - p_0 q_0  - \int_0^h \left[ p(t) \dot{q}(t)-H(q(t), p(t)) \right] dt \label{exact_LeftHd},
\end{align}
where $(q,p)$ satisfies Hamilton's equations with boundary conditions $q(h)=q_h$ and $p(0)=p_0$. A Type III Hamiltonian variational integrator is constructed by using an approximate discrete left Hamiltonian $H_d^-$, and applying the \textbf{discrete left Hamilton's equations},
\begin{equation} \label{Discrete Left Eq}
	p_1= - D_1H_d^-(q_1,p_0), \qquad q_0 = - D_2H_d^-(q_1,p_0),
\end{equation}
which implicitly define the integrator, $\tilde{F}_{H_d^-}:(q_0,p_0) \mapsto (q_1,p_1)$.

We now derive analogous results to those of Section~\ref{section: Constrained Variational Lagrangian} from the Hamiltonian perspective. As in the Lagrangian case, we will assume we have a configuration manifold $\mathcal{M}$, a holonomic constraint function $\mathcal{C} : \mathcal{M} \rightarrow \mathbb{R}^d$, and that the dynamics are constrained to the submanifold $\mathcal{Q} = \mathcal{C}^{-1}(0)$.

\subsection{Continuous Constrained Variational Hamiltonian Mechanics}

The following theorem presents the equivalence between a continuous constrained variational principle and continuous constrained Hamilton's equations in the Type II case, generalizing Lemma 2.1 from~\cite{LeZh2011} to include holonomic constraints:

\begin{theorem} \label{theorem: Action and Equations II}
	Consider the Type II \textbf{constrained action functional} $\mathfrak{S} : C^2([0,T],T^*\mathcal{Q} \times \Lambda) \rightarrow \mathbb{R}$
	\begin{equation}\label{Constrained Action Functional}
		\mathfrak{S} (q(\cdot),p(\cdot), \lambda(\cdot)) = p(T)q(T) - \int_{0}^{T}{ \left[ p(t) \dot{q}(t) - H(q(t),p(t)) - \langle \lambda(t) , \mathcal{C}(q(t)) \rangle \right] dt}.
	\end{equation}
	The condition that $\mathfrak{S} (q(\cdot),p(\cdot), \lambda(\cdot)) $ is stationary with respect to the boundary conditions $\delta q(0) = 0$ and $\delta p(T) =0$ is equivalent to $(q(\cdot),p(\cdot),\lambda(\cdot))$ satisfying \textbf{Hamilton's constrained equations}
	\begin{equation} \label{HamiltonConstEqs}
		\dot{q} = \frac{\partial H}{\partial p} (q,p) , \qquad \dot{p} = -\frac{\partial H}{\partial q} (q,p) - \langle \lambda , \nabla  \mathcal{C} (q) \rangle  ,  \qquad \mathcal{C}(q) = 0 .
	\end{equation}
\proof{See Appendix~\ref{Appendix: Type II Action}.  \qed}
\end{theorem}

As in the Type II case, we can derive a theorem relating a continuous constrained variational principle and continuous constrained Hamilton's equations in the Type III case:

\begin{theorem} \label{theorem: Action and Equations III}
	Consider the Type III \textbf{constrained action functional} $\mathfrak{S} : C^2([0,T],T^*\mathcal{Q} \times \Lambda) \rightarrow \mathbb{R}$
	\begin{equation}\label{Constrained Action Functional III}
		\mathfrak{S} (q(\cdot),p(\cdot), \lambda(\cdot)) = - p(0)q(0) - \int_{0}^{T}{ \left[ p(t) \dot{q}(t) - H(q(t),p(t)) - \langle \lambda(t) , \mathcal{C}(q(t)) \rangle \right] dt}.
	\end{equation}
	The condition that $\mathfrak{S} (q(\cdot),p(\cdot), \lambda(\cdot)) $ is stationary with respect to the boundary conditions $\delta q(T) = 0$ and $\delta p(0) =0$ is equivalent to $(q(\cdot),p(\cdot),\lambda(\cdot))$ satisfying \textbf{Hamilton's constrained equations}
	\begin{equation} \label{HamiltonConstEqsIII}
		\dot{q} = \frac{\partial H}{\partial p} (q,p) , \qquad \dot{p} = -\frac{\partial H}{\partial q} (q,p) - \langle \lambda , \nabla  \mathcal{C} (q) \rangle,  \qquad \mathcal{C}(q) = 0 .
	\end{equation}
	\proof{See Appendix~\ref{Appendix: Type III Action}.  \qed}
\end{theorem}

\begin{remark}
	Hamilton's constrained equations are the same in the Type II and Type III formulations of Hamiltonian mechanics, and they can be thought of as the Hamilton's equations generated by the augmented Hamiltonian 
	\begin{equation}
		\bar{H} \left( q,\lambda, p, \mathcal{p} \right)  = H(q,p ) + \langle \lambda , \mathcal{C} (q) \rangle,
	\end{equation} where $\mathcal{p}$ is the conjugate momentum for the variable $\lambda$. Furthermore, they are equivalent to the constrained Euler--Lagrange equations \eqref{ELConstEqs}, provided that the Lagrangian $L$ is hyperregular.
\end{remark}

\begin{remark}
	It is sometimes beneficial to augment the continuous equations with the equation $ \langle \frac{\partial H}{\partial p} (q,p)  , \nabla \mathcal{C}(q)  \rangle  =0$ (and analogously for the discrete case) to ensure that the momentum $p$ lies in the cotangent space to the manifold, as explained and illustrated in~\cite[Chapter VII]{HaLuWa2006}. 
\end{remark}

We will now generalize Theorem 2.2 and its Type III analogue from~\cite{LeZh2011} to include holonomic constraints $\mathcal{C}(q)$ using Lagrange multipliers $\lambda : [0,T] \rightarrow \Lambda$.

In the Type II case, consider the function $\mathcal{S}(q_0,p_T)$ given by the extremal value of the constrained action functional $\mathfrak{S}$ over the family of curves $(q(\cdot),p(\cdot) , \lambda(\cdot)) $ satisfying the boundary conditions $q(0)=q_0$ and $p(T)=p_T$:
\begin{equation}
	\mathcal{S}(q_0,p_T) = \ext_{ \substack{(q,p,\lambda)\in C^2([0,T],T^*\mathcal{Q} \times \Lambda)   \\ q(0)=q_0, \quad  p(T)=p_T}}{\mathfrak{S} (q(\cdot),p(\cdot) , \lambda(\cdot)) }. 
\end{equation}
The following theorem shows that $\mathcal{S}(q_0,p_T)$ is a generating function for the flow of the continuous constrained Hamilton's equations:

\begin{theorem} \label{theorem: Flow Map and Equations II}
	The exact time-$T$ flow map of Hamilton's equations $(q_0,p_0) \mapsto (q_T,p_T)$ is implicitly given by the following relations:
	\begin{equation}
		q_T = D_2 \mathcal{S} (q_0,p_T), \qquad p_0 = D_1 \mathcal{S} (q_0,p_T).
	\end{equation}
	In particular, $ \mathcal{S} (q_0,p_T)$ is a Type II generating function that generates the exact flow of the constrained Hamilton's equations \eqref{HamiltonConstEqs}. 
	
	\proof{See Appendix~\ref{Appendix: Type II Flow}. \qed
}
\end{theorem}

In the Type III case, consider the function $\mathcal{S}(q_T,p_0)$ given by the extremal value of the constrained action functional $\mathfrak{S}$ over the family of curves $(q(\cdot),p(\cdot) , \lambda(\cdot)) $ satisfying the boundary conditions $q(T)=q_T$ and $p(0)=p_0$:
\begin{equation}
	\mathcal{S}(q_T,p_0) = \ext_{ \substack{(q,p,\lambda)\in C^2([0,T],T^*\mathcal{Q} \times \Lambda)   \\ q(T)=q_T, \quad  p(0)=p_0}}{\mathfrak{S} (q(\cdot),p(\cdot) , \lambda(\cdot)) }.
\end{equation}
The following theorem shows that $\mathcal{S}(q_T,p_0)$ is a generating function for the flow of the continuous constrained Hamilton's equations:

\begin{theorem}  \label{theorem: Flow Map and Equations III}
	The exact time-$T$ flow map of Hamilton's equations $(q_0,p_0) \mapsto (q_T,p_T)$ is implicitly given by the following relations:
	\begin{equation}
		q_0 = -  D_2 \mathcal{S} (q_T,p_0), \qquad p_T = - D_1 \mathcal{S} (q_T,p_0).
	\end{equation}
	In particular, $ \mathcal{S} (q_T,p_0)$ is a Type III generating function that generates the exact flow of the constrained Hamilton's equations \eqref{HamiltonConstEqsIII}. 
	\proof{See Appendix~\ref{Appendix: Type III Flow}. \qed \\
	}
\end{theorem}

\subsection{Discrete Constrained Variational Hamiltonian Mechanics}

Let us now extend the results of Section 3 from~\cite{LeZh2011} to introduce a discrete formulation of variational Hamiltonian mechanics which includes holonomic constraints. Suppose we are given a partition $0 = t_0 < t_1 < \ldots < t_N = T$ of the interval $[0,T]$, and a discrete curve in $T^* \mathcal{Q} \times \Lambda$, denoted by $\{  (q_k,p_k, \lambda _k) \}_{k=0}^{N}$, such that $q_k \approx q(t_k)$, $p_k \approx p(t_k)$ and $\lambda_k \approx \lambda(t_k)$.

We formulate discrete constrained variational Hamiltonian mechanics in terms of the following discrete analogues of the constrained action functional $\mathfrak{S}$ given by equation \eqref{Constrained Action Functional}:
	\begin{align} \label{Discrete Constrained Action Functional}
	\mathfrak{S}_d^+ \left(\{  (q_k,p_k, \lambda _k) \}_{k=0}^{N} \right) & = p_N q_N - \sum_{k=0}^{N-1}{\left[ p_{k+1} q_{k+1} - H_d^+(q_k,p_{k+1}) - \langle \lambda_k , \mathcal{C}(q_k) \rangle \right]}, 	\\
	\mathfrak{S}_d^- \left(\{  (q_k,p_k, \lambda _k) \}_{k=0}^{N} \right) & = -p_0 q_0 - \sum_{k=0}^{N-1}{\left[ - p_{k} q_{k} - H_d^-(q_{k+1},p_k) - \langle \lambda_{k+1} , \mathcal{C} (q_{k+1}) \rangle \right]},
\end{align}
where
\begin{align}\label{Discrete Right Hamiltonian}
	H_d^+(q_k,p_{k+1}) & \approx  \ext_{ \substack{(q,p,\lambda)\in C^2([t_k,t_{k+1}],T^*\mathcal{Q} \times \Lambda)   \\ q(t_k)=q_k, \quad  p(t_{k+1})=p_{k+1}}} p(t_{k+1}) q(t_{k+1})   - \int_{t_k}^{t_{k+1}}{\left[   p(t) \dot{q}(t) - H(q(t),p(t))  \right] dt} \\ \label{Discrete Left Hamiltonian}
	H_d^-(q_{k+1},p_k) & \approx  \ext_{ \substack{(q,p,\lambda)\in C^2([t_k,t_{k+1}],T^*\mathcal{Q} \times \Lambda)   \\ q(t_{k+1})=q_{k+1}, \quad  p(t_{k})=p_{k}}}   -p(t_{k}) q(t_{k})    - \int_{t_k}^{t_{k+1}}{\left[   p(t) \dot{q}(t) - H(q(t),p(t))  \right] dt} . 
\end{align}

We can now derive discrete analogues of Theorems~\ref{theorem: Action and Equations II} and~\ref{theorem: Action and Equations III} relating discrete variational principles to discrete constrained Hamilton's equations, generalizing Lemma 3.1 from~\cite{LeZh2011}:

	\begin{theorem}  \label{theorem: Discrete Variational Principle II}
		The Type II discrete Hamilton's phase space variational principle
		\begin{equation}
			\delta \mathfrak{S}_d^{+} \left(\{  (q_k,p_k, \lambda _k) \}_{k=0}^{N} \right)  = 0
		\end{equation}
		is equivalent to the \textbf{discrete constrained right Hamilton's equations}
		\begin{equation} \label{eq: discrete constrained right Hamilton equations}
			q_{k+1}  = D_2 H_d^+(q_{k},p_{k+1}),  \qquad p_k  =  D_1H_d^+(q_k,p_{k+1}) + \langle \lambda_k, \nabla\mathcal{C}(q_k)  \rangle, \qquad  \mathcal{C}(q_k)  =0, 
		\end{equation}
		where $H_d^+ (q_k,p_{k+1}) $ is defined via equation \eqref{Discrete Right Hamiltonian}.
		\proof{See Appendix~\ref{Appendix: Type II Discrete Variational}. \qed
		}
	\end{theorem}

	\begin{theorem}  \label{theorem: Discrete Variational Principle III}
	The Type III discrete Hamilton's phase space variational principle
	\begin{equation}
		\delta \mathfrak{S}_d^{-} \left(\{  (q_k,p_k, \lambda _k) \}_{k=0}^{N} \right)  = 0
	\end{equation}
	is equivalent to the \textbf{discrete constrained left Hamilton's equations}
	\begin{equation}
		q_{k}  = -D_2 H_d^-(q_{k+1},p_{k}),  \qquad p_{k+1}  =  -D_1H_d^-(q_{k+1},p_{k}) - \langle \lambda_{k+1}, \nabla\mathcal{C}(q_{k+1})  \rangle, \qquad  \mathcal{C}(q_k)  =0, 
	\end{equation}
	where $H_d^- (q_{k+1},p_{k}) $ is defined via equation \eqref{Discrete Left Hamiltonian}.
	\proof{See Appendix~\ref{Appendix: Type III Discrete Variational}. \qed
	}
\end{theorem}

\begin{remark}
	These discrete constrained Hamilton's equations can be thought of as the discrete Hamilton's equations generated by the augmented discrete Hamiltonians
	\begin{align} 
			\bar{H}_d^{+} \left((q_k,\lambda_k),(p_{k+1},\mathcal{p}_{k+1})\right) & = H_d^{+} (q_k,p_{k+1}) + \langle \lambda_k , \mathcal{C}(q_k) \rangle, \\
		\bar{H}_d^{-}\left((q_{k+1},\lambda_{k+1}),(p_{k},\mathcal{p}_{k})\right) &= H_d^{-} (q_k,p_{k+1}) + \langle \lambda_{k+1} , \mathcal{C}(q_{k+1}) \rangle. 
	\end{align}
This augmented Hamiltonian perspective together with the augmented Lagrangian perspective from Remark~\ref{remark: augmented Lagrangian} imply that the constrained $\bar{H}_d^+$ variational integrator is equivalent to the constrained $\bar{L}_d^+$ variational integrator whenever the $H_d^+$ variational integrator is equivalent to the $L_d^+$ variational integrator (and similarly for the integrators generated by $\bar{H}_d^-$ and $\bar{L}_d^-$). Examples where this happens are presented in~\cite{ScShLe2017} for Taylor variational integrators provided the Lagrangian is hyperregular, and in~\cite{LeZh2011} for generalized Galerkin variational integrators constructed using the same choices of basis functions and numerical quadrature formula provided the  Hamiltonian is hyperregular. \\
\end{remark}

\section{Error Analysis for Variational Integrators} \label{sec: Error Analysis}
\subsection{Unconstrained Error Analysis}

Theorem 2.3.1 of~\cite{MaWe2001} states that if a discrete Lagrangian, $L_d:\mathcal{Q}\times \mathcal{Q}\rightarrow\mathbb{R}$, approximates the exact discrete Lagrangian $L_d^E:\mathcal{Q}\times \mathcal{Q}\rightarrow\mathbb{R}$ to order $r$, i.e.,
\begin{align} 
	L_d(q_0, q_h)=L_d^E(q_0,q_h)+\mathcal{O}(h^{r+1}) ,
\end{align}
then the discrete Hamiltonian map $\tilde{F}_{L_d}:(q_k,p_k)\mapsto(q_{k+1},p_{k+1})$, viewed as a one-step method defined implicitly from the discrete Euler--Lagrange equations
\begin{equation} \label{eq: DEL ErrorAnalysis}
	D_1 L_d(q_k, q_{k+1}) +D_2 L_d(q_{k-1}, q_{k}) =0, 
\end{equation}
or equivalently in terms of the implicit discrete Euler--Lagrange equations, which involve the corresponding discrete momenta via the discrete Legendre transforms,
\begin{equation}
	p_k= -D_1L_d(q_k,q_{k+1}), \qquad p_{k+1}=D_2L_d(q_{k},q_{k+1}),
\end{equation}
has order of accuracy $r$.

Theorem 2.3.1 of~\cite{MaWe2001} has an analogue for Hamiltonian variational integrators. Theorem 2.2 in~\cite{ScLe2017} states that if a discrete right Hamiltonian $H^+_d$ approximates the exact discrete right Hamiltonian $H_d^{+,E}$ to order $r$, i.e.,
\begin{align} 
	H^+_d(q_0, p_h)=H_d^{+,E}(q_0,p_h)+\mathcal{O}(h^{r+1}),
\end{align}
then the discrete right Hamiltonian map $\tilde{F}_{H^+_d}:(q_k,p_k)\mapsto(q_{k+1},p_{k+1})$, viewed as a one-step method defined implicitly by the discrete right Hamilton's equations
\begin{equation}
	p_k=D_1H_d^+(q_k,p_{k+1}), \qquad q_{k+1}=D_2H_d^+(q_k,p_{k+1}),
\end{equation}
is order $r$ accurate. As mentioned in~\cite{ScLe2017}, the proof of Theorem 2.2 in~\cite{ScLe2017} can easily be adjusted to prove an equivalent theorem for the discrete left Hamiltonian case, which states that if a discrete left Hamiltonian $H^-_d$ approximates the exact discrete left Hamiltonian $H_d^{-,E}$ to order $r$, i.e.,
\begin{align} 
	H^-_d(q_1, p_0)=H_d^{-,E}(q_1,p_0)+\mathcal{O}(h^{r+1}),
\end{align}
then the discrete left Hamiltonian map $\tilde{F}_{H^-_d}:(q_k,p_k)\mapsto(q_{k+1},p_{k+1})$, viewed as a one-step method defined implicitly by the discrete left Hamilton's equations
\begin{equation} 
	p_{k+1}= - D_1H_d^-(q_{k+1},p_k), \qquad q_k = - D_2H_d^-(q_{k+1},p_k),
\end{equation}
is order $r$ accurate. Many other properties of the integrator, such as momentum conservation properties of the method, can be determined by analyzing the associated discrete Lagrangian or Hamiltonian, as opposed to analyzing the integrator directly.  We will exploit these error analysis results to derive analogous results for the constrained versions discussed in Sections~\ref{section: Constrained Variational Lagrangian} and~\ref{Sec: Constrained H Mechanics}.

\subsection{Constrained Error Analysis} \label{sec: Constrained Error Analysis}  

For the Lagrangian case, we can think of the Lagrange multipliers $\lambda$ as extra position coordinates and define an augmented Lagrangian $\bar{L}$ via
\begin{equation} \label{eq: Augmented Lagrangian}
	\bar{L}\left((q,\lambda),(\dot{q}, \dot{\lambda})\right)  =  L(q,\dot{q}) - \langle \lambda , \mathcal{C}(q) \rangle.
\end{equation}

\noindent A corresponding augmented discrete Lagrangian is given by
\begin{equation}
	\bar{L}_d \left((q_k,\lambda_k),(q_{k+1},\lambda_{k+1})\right) = L_d(q_k,q_{k+1}) - \langle \lambda_k , \mathcal{C}(q_k) \rangle ,
\end{equation}
and the discrete Euler--Lagrange equations \eqref{eq: DEL ErrorAnalysis}
\begin{equation} 
	D_1 \bar{L}_d \left((q_k,\lambda_k),(q_{k+1},\lambda_{k+1})\right) +D_2 \bar{L}_d \left((q_{k-1},\lambda_{k-1}),(q_{k},\lambda_{k})\right) =0, 
\end{equation}
yield the discrete constrained Euler--Lagrange equations
\begin{equation}
	 D_1 L_d(q_k,q_{k+1})  + D_2 L_d(q_{k-1},q_{k}) = \langle \lambda_k ,\nabla \mathcal{C}(q_k) \rangle,  \qquad \mathcal{C}(q_k)  =0,  
\end{equation}
derived in Section~\ref{sec: Discrete Lagrangian}. As a consequence, we can apply Theorem 2.3.1 of~\cite{MaWe2001} to the augmented Lagrangian \eqref{eq: Augmented Lagrangian} and obtain the following result:
\begin{theorem}
	Suppose that for an exact discrete Lagrangian $L_d^E$ and a discrete Lagrangian $L_d$, 
	\begin{align} 
		L_d(q_0, q_h)  - \langle \lambda_0 , \mathcal{C}(q_0) \rangle  =L_d^E(q_0,q_h) -  \int_{0}^{h}{\langle \lambda(t) , \mathcal{C}(q(t)) \rangle dt} +\mathcal{O}(h^{r+1}).
	\end{align}
Then, the discrete map $(q_k,p_k,\lambda_k)\mapsto(q_{k+1},p_{k+1},\lambda_{k+1})$, viewed as a one-step method defined implicitly by the discrete constrained Euler--Lagrange equations, has order of accuracy $r$.
\end{theorem}

For the Hamiltonian case, we can think of the Lagrange multipliers $\lambda$ as extra position coordinates and define conjugate momenta $\mathcal{p}$, which are constants of motion since the time-derivative of $\lambda$ does not appear anywhere, and are constrained to be zero. The augmented Hamiltonian $\bar{H}$, given by
\begin{equation} 
	\bar{H}\left((q,\lambda),(p,\mathcal{p} )\right)  =  H(q,p) + \langle \lambda , \mathcal{C}(q) \rangle,
\end{equation}
yields the following augmented left and right discrete Hamiltonians
\begin{align}
	\bar{H}_d^{-} \left((q_{k+1},\lambda_{k+1}),(p_{k},\mathcal{p}_{k})\right) &= H_d^{-} (q_{k+1},p_k) + \langle \lambda_{k+1} , \mathcal{C}(q_{k+1}) \rangle ,\\
	\bar{H}_d^{+} \left((q_k,\lambda_k),(p_{k+1},\mathcal{p}_{k+1})\right) &= H_d^{+} (q_k,p_{k+1}) + \langle \lambda_k , \mathcal{C}(q_k) \rangle ,
\end{align}
and the discrete left and right Hamilton's equations
\begin{alignat}{2} 
	(p_{k+1},\mathcal{p}_{k+1})&= - D_1\bar{H}_d^-\left((q_{k+1},\lambda_{k+1}),(p_k,\mathcal{p}_k)\right),& \,\, (q_k,\lambda_k) &= - D_2\bar{H}_d^-\left((q_{k+1},\lambda_{k+1}),(p_k,\mathcal{p}_k)\right),\\
(p_k,\mathcal{p}_k)&=D_1\bar{H}_d^+\left((q_k,\lambda_k),(p_{k+1},\mathcal{p}_{k+1})\right),& \,\, (q_{k+1},\lambda_{k+1})&=D_2\bar{H}_d^+\left((q_k,\lambda_k),(p_{k+1},\mathcal{p}_{k+1})\right),
\end{alignat}
yield the discrete constrained left Hamilton's equations
\begin{equation}
	q_k  = - D_2 H_d^-(q_{k+1},p_{k}),  \qquad p_{k+1}  =  - D_1H_d^-(q_{k+1},p_{k}) - \langle \lambda_{k+1} ,\nabla \mathcal{C}(q_{k+1})  \rangle, \quad \mathcal{C}(q_k)  =0, 
\end{equation}
and the discrete constrained right Hamilton's equations
\begin{equation}
	q_{k+1}  = D_2 H_d^+(q_{k},p_{k+1}),  \qquad p_k  =  D_1H_d^+(q_k,p_{k+1}) +  \langle \lambda_k, \nabla  \mathcal{C}(q_k)  \rangle, \quad  \mathcal{C}(q_{k})  =0,
\end{equation}
derived in Section~\ref{Sec: Constrained H Mechanics}. As a consequence, we can apply Theorem 2.2 in~\cite{ScLe2017} and its Type III analogue to the augmented Hamiltonians and obtain the following results
\begin{theorem} \label{Thm: Type II Error Analysis}
	Suppose that given an exact discrete right Hamiltonian $H_d^{+,E}$ and a discrete right Hamiltonian $H_d^+$, we have
	\begin{align} 
			H^+_d(q_0, p_h) +  \langle \lambda_0 , \mathcal{C}(q_0) \rangle  =H_d^{+,E}(q_0,p_h)+  \int_{0}^{h}{\langle \lambda(t) , \mathcal{C}(q(t)) \rangle dt} +\mathcal{O}(h^{r+1}).
	\end{align}
	Then, the discrete map $(q_k,p_k,\lambda_k)\mapsto(q_{k+1},p_{k+1},\lambda_{k+1})$, viewed as a one-step method defined implicitly by the discrete constrained right Hamilton's equations, has order of accuracy $r$.
\end{theorem}
\begin{theorem}
	Suppose that given an exact discrete left Hamiltonian $H_d^{-,E}$ and a discrete left Hamiltonian $H_d^-$, we have
	\begin{align} 
		H^-_d(q_h, p_0) +  \langle \lambda_h , \mathcal{C}(q_h) \rangle  =H_d^{-,E}(q_h,p_0)+  \int_{0}^{h}{\langle \lambda(t) , \mathcal{C}(q(t)) \rangle dt} +\mathcal{O}(h^{r+1}).
	\end{align}
	Then, the discrete map $(q_k,p_k,\lambda_k)\mapsto(q_{k+1},p_{k+1},\lambda_{k+1})$, viewed as a one-step method defined implicitly by the discrete constrained left Hamilton's equations, has order of accuracy $r$. \\
\end{theorem} 

\section{Variational Riemannian Accelerated Optimization} \label{sec: Optimization Section}

\subsection{Riemannian Geometry}

We first introduce the main notions from Riemannian geometry that will be used throughout this section (see~\cite{Absil2008,Boumal2020,Duruisseaux2021Riemannian,Jost2017,Lee2019,Lang1999} for more details). 

\begin{definition} \label{def: fiber metric}
	Suppose we have a Riemannian manifold $\mathcal{Q}$ with Riemannian metric  $g(\cdot,\cdot) = \langle \cdot , \cdot \rangle$, represented by the positive-definite symmetric matrix  $(g_{ij}) $ in local coordinates. Then, we define the \textbf{musical isomorphism} $g^{\flat}:T\mathcal{Q} \rightarrow T^*\mathcal{Q}$ via 
	$$ g^{\flat}(u)(v) = g_q(u,v)  \quad  \forall  q\in \mathcal{Q}  \text{ and }  \forall u,v\in T_q\mathcal{Q},  $$ and its \textbf{inverse musical isomorphism} $g^{\sharp}:T^*\mathcal{Q} \rightarrow T\mathcal{Q}$. The Riemannian metric $g(\cdot,\cdot) = \langle \cdot , \cdot \rangle$ induces a \textbf{fiber metric} $g^*(\cdot ,\cdot) =   \llangle \cdot , \cdot \rrangle $ on $T^* \mathcal{Q}$ via
	$$  \llangle u , v \rrangle = \langle g^{\sharp}(u), g^{\sharp}(v) \rangle  \quad \forall u,v \in T^* \mathcal{Q}, $$
	represented by the positive-definite symmetric matrix $(g^{ij})$ in local coordinates, which is the inverse of the Riemannian metric matrix $(g_{ij}) $. 
\end{definition}

\begin{definition}
	Denoting the differential of $f$ by $df$, the \textbf{Riemannian gradient} $\emph{grad}f(q) \in T_q  \mathcal{Q}$ at a point $q\in \mathcal{Q}$ of a smooth function $f:\mathcal{Q} \rightarrow \mathbb{R}$ is the tangent vector at $q$ such that 
	$$ \langle  \emph{grad}f(q) , u \rangle = df(q) u   \qquad \forall u\in T_q \mathcal{Q}. $$
	This can also be expressed in terms of the inverse musical isomorphism, $\emph{grad}f(q)=g^\sharp(df(q))$.
\end{definition}

\begin{definition}
	A \textbf{geodesic} in a Riemannian manifold $\mathcal{Q}$ is a parametrized curve $\gamma : [0,1] \rightarrow \mathcal{Q}$ which is of minimal local length, and is a generalization of the notion of straight line from Euclidean spaces to Riemannian manifolds. The other generalization of straight lines involves curves having zero ``acceleration" or constant ``speed", which requires the introduction of an affine connection. These two generalizations are equivalent if the Riemannian manifold is endowed with the Levi--Civita connection. Given two points $q,\tilde{q} \in \mathcal{Q}$, a vector in $T_q \mathcal{Q}$ can be transported to $T_{\tilde{q}}\mathcal{Q}$ along a geodesic $\gamma$ by an operation $\Gamma_q^{\tilde{q}}:T_q \mathcal{Q}\rightarrow T_{\tilde{q}}\mathcal{Q}$ called the \textbf{parallel transport along} $\gamma$. 
\end{definition}

\begin{definition}
	The \textbf{Riemannian Exponential map} $\emph{Exp}_q:T_q \mathcal{Q} \rightarrow \mathcal{Q}$ at $q\in \mathcal{Q}$ is defined via
	$$ \emph{Exp}_q(v) = \gamma_v(1),  $$
	where $\gamma_v$ is the unique geodesic in $\mathcal{Q}$ such that $\gamma_v(0) = q$ and $\gamma_v'(0) = v$, for any $v\in T_q \mathcal{Q} $.  $\emph{Exp}_q$ is a diffeomorphism in some neighborhood $ U \subset T_q\mathcal{Q}$ containing 0, so we can define its inverse map, the \textbf{Riemannian Logarithm map} $\emph{Log}_p : \emph{Exp}_q(U) \rightarrow T_q \mathcal{Q}$.
\end{definition}

\begin{definition}
	A \textbf{retraction} on a manifold $\mathcal{Q}$ is a smooth mapping $\mathcal{R}: T\mathcal{Q} \rightarrow \mathcal{Q}$ such that for any $q \in \mathcal{Q}$, the restriction $\mathcal{R}_q : T_q\mathcal{Q} \rightarrow \mathcal{Q} $ of $\mathcal{R}$ to $T_q\mathcal{Q} $ satisfies
	\begin{itemize}
		\item $\mathcal{R}_q(0_q) = q$, where $0_q$ denotes the zero element of $T_q\mathcal{Q} $, 
		\item $T_{0_q}\mathcal{R}_q = \mathbb{I}_{T_q\mathcal{Q} }$ with the canonical identification  $T_{0_q}T_{q}\mathcal{Q} \simeq T_{q}\mathcal{Q}$, where $T_{0_q}\mathcal{R}_q$ is the tangent map of $\mathcal{R}$ at $0_q \in T_{q}\mathcal{Q}$ and $\mathbb{I}_{T_q\mathcal{Q} }$ is the identity map on $T_{q}\mathcal{Q}$.
	\end{itemize}
			The Riemannian Exponential map is a natural example of a retraction on a Riemannian manifold.
\end{definition}

\begin{definition}
	A subset $A$ of a Riemannian manifold $\mathcal{Q}$ is called \textbf{geodesically uniquely convex} if every two points of $A$ are connected by a unique geodesic in $A$. A function $f:\mathcal{Q} \rightarrow \mathbb{R}$ is called \textbf{geodesically convex} if for any two points $q,\tilde{q} \in \mathcal{Q}$ and a geodesic $\gamma$ connecting them, 
	$$ f(\gamma(t)) \leq (1-t) f(q) +t f(\tilde{q})  \qquad \forall t\in [0,1].$$ Note that if $f$ is a smooth geodesically convex function on a geodesically uniquely convex subset $A$,
	$$ f(q) - f(\tilde{q}) \geq \langle  \emph{grad}f(\tilde{q}) , \emph{Log}_{\tilde{q}}(q) \rangle   \qquad \forall q,\tilde{q} \in A.$$
	A function $f:A\rightarrow \mathbb{R}$ is called \textbf{geodesically $\alpha$-weakly-quasi-convex} ($\alpha$-WQC) with respect to $q \in \mathcal{Q}$ for some $\alpha \in (0,1]$ if
	$$ \alpha \left(f(q) - f(\tilde{q})\right) \geq \langle  \emph{grad}f(\tilde{q}) , \emph{Log}_{\tilde{q}}(q) \rangle   \qquad \forall \tilde{q} \in A.$$
	Note that a local minimum of a geodesically convex or $\alpha$-WQC function is also a global minimum.
\end{definition}

\begin{definition}
	Given a Riemannian manifold $\mathcal{Q}$ with sectional curvature bounded below by $K_{\min}$, and an upper bound $D$ for the diameter of the domain of consideration, define
	\begin{align}\label{eq: zeta}
		\zeta = 
		\begin{cases}
			\sqrt{-K_{\min}} D \coth{ (\sqrt{-K_{\min}} D) }  & \quad \text{if }  K_{\min} < 0 \\ 1   & \quad \text{if }  K_{\min} \geq 0 
		\end{cases} .
	\end{align}
	Note that $\zeta \geq 1$ since $x \coth{x} \geq 1$ for all real values of $x$. 
\end{definition}

\subsection{Hamiltonian Approach}

Our approach consists in integrating the Riemannian Bregman Hamiltonian systems derived in~\cite{Duruisseaux2021Riemannian} which evolve on the Riemannian manifold $\mathcal{Q}$, via discrete constrained variational Hamiltonian integrators which enforce the numerical solution to lie on the Riemannian manifold $\mathcal{Q}$. With $\zeta$ given by equation \eqref{eq: zeta}, we know from~\cite{Duruisseaux2021Riemannian} that if we let $\lambda = \zeta$ in the geodesically convex case, and $\lambda = \zeta / \alpha$ in the geodesically $\alpha$-weakly-quasi-convex case, we obtain the Direct approach Riemannian $p$-Bregman Hamiltonian
\begin{equation}  \label{H-R-Direct}
	\bar{\mathcal{H}}_{p}(\bar{Q},\bar{R})  = 	 \frac{p}{2(Q^t)^{\lambda p +1}} \llangle R , R\rrangle  + Cp(Q^t)^{(\lambda +1)p-1} f(Q) + R^t,
\end{equation}
and the Adaptive approach Riemannian $p\rightarrow \mathring{p}$ Bregman Hamiltonian
\begin{equation}  \label{H-R-Adaptive}
	\bar{\mathcal{H}}_{p \rightarrow \mathring{p}}(\bar{Q},\bar{R})  = 	 \frac{p^2}{2\mathring{p} (Q^t)^{\lambda p +\mathring{p}/p}}  \llangle R , R\rrangle  + \frac{Cp^2}{\mathring{p}}(Q^t)^{(\lambda +1)p-\mathring{p}/p} f(Q)+ \frac{p}{\mathring{p}} (Q^t)^{1-\mathring{p}/p}   R^t.
\end{equation}
It was proven in~\cite{Duruisseaux2021Riemannian} that along the trajectories of the Riemannian $p$-Bregman dynamics, $f(Q(t))$ converges to its optimal value at a rate of $\mathcal{O}(1/t^p)$, under suitable assumptions on $\mathcal{Q}$.

\begin{remark}
In the vector space setting, these Riemannian Bregman Hamiltonians reduce to the Direct and Adaptive approach Bregman Hamiltonians derived in~\cite{duruisseaux2020adaptive} for convex functions:
	\begin{equation}\label{H-Direct}
		\bar{H}_{p}(\bar{q},\bar{r}) = \frac{p}{2(q^t)^{p+1}} \langle r , r \rangle + Cp(q^t)^{2p-1} f(q) + r^t ,
	\end{equation}
	\begin{equation} \label{H-Adaptive}
		\bar{H}_{p\rightarrow \mathring{p}}(\bar{q},\bar{r}) =   \frac{p^2}{2\mathring{p}(q^t)^{p+\mathring{p}/p} }  \langle r , r \rangle  + \frac{Cp^2}{\mathring{p}} (q^t)^{2p-\mathring{p}/p} f(q) + \frac{p}{\mathring{p}} (q^t)^{1-\mathring{p}/p} r^t  . 
	\end{equation}
\end{remark}

\hfill 

\subsection{Some Optimization Problems on Riemannian Manifolds}  \label{sec: List of Problems}

\subsubsection{Rayleigh Quotient Minimization on the Unit Sphere}

An eigenvector $v$ corresponding to the largest eigenvalue of a symmetric $n \times n$ matrix $A$ maximizes the Rayleigh quotient $\frac{v^\top Av}{ v^\top v}$ over $\mathbb{R}^n$. Thus, a unit eigenvector corresponding to the largest eigenvalue of the matrix $A$ is a minimizer of the function $f(v) = -  v^\top Av$ over the unit sphere $\mathcal{Q} = \mathbb{S}^{n-1}$, which can be thought of as a Riemannian submanifold with constant positive curvature $K=1$ of $\mathbb{R}^{n}$ endowed with the Riemannian metric inherited from the Euclidean inner product $g_v(u,w) = u^\top w$. Solving the Rayleigh quotient optimization problem efficiently is challenging when the given symmetric matrix $A$ is ill-conditioned and high-dimensional. Note that an efficient algorithm that solves the above minimization problem can also be used to find eigenvectors corresponding to the smallest eigenvalue of $A$ by using the fact that the eigenvalues of $A$ are the negative of the eigenvalues of $-A$.

\subsubsection{Eigenvalue and Procrustes Problems on the Stiefel Manifold}

When endowed with the Riemannian metric $g_X(A,B) = \text{Trace}(A^\top B)$, the \textbf{Stiefel manifold}
\begin{equation} \text{St}(m,n) = \{X\in \mathbb{R}^{n\times m} | X^\top X= I_m \} 
\end{equation}
is a Riemannian submanifold of $\mathbb{R}^{n\times m}$. The tangent space at any $X \in \text{St}(m,n)$ is given by $T_X \text{St}(m,n) = \{ Z\in \mathbb{R}^{n\times m} | X^\top Z+Z^\top X=0  \},$
and the orthogonal projection $P_X$ onto $T_X \text{St}(m,n)$ is given by $P_X Z = Z - \frac{1}{2} X(X^\top Z + Z^\top X).$ A retraction on $\text{St}(m,n) $ is given by $\mathcal{R}_X(\xi) = \text{qf} (X+\xi) ,$
where $\text{qf}(A)$ denotes the $Q$ factor of the QR factorization of the matrix $A\in \mathbb{R}^{n\times m}$ as $A=QR$ where $Q\in \text{St}(m,n)$ and $R$ is an upper triangular $n\times m$ matrix with strictly positive diagonal elements~\cite{Absil2008}.

A generalized eigenvector problem consists of finding the $m$ smallest eigenvalues of a $n\times n$ symmetric matrix $A$ and corresponding eigenvectors. This problem can be formulated as a Riemannian optimization problem on the Stiefel manifold $\text{St}(m,n)$ via the Brockett cost function
\begin{equation}
	f:\text{St}(m,n) \rightarrow \mathbb{R}, \quad  X\mapsto f(X) = \text{Trace}(X^\top AXN),
\end{equation}
where $N = \text{diag}(\mu_1 , \ldots , \mu_m)$ for arbitrary $0 \leq \mu_1 \leq \ldots \leq \mu_m $. The columns of a global minimizer of $f$ are eigenvectors corresponding to the $m$ smallest eigenvalues of $A$ (see~\cite{Absil2008}).  If we define $\bar{f} : \mathbb{R}^{n\times m} \rightarrow \mathbb{R}$ via $X\mapsto \bar{f}(X) = \text{Trace}(X^\top AXN),$ then $f$ is the restriction of $\bar{f}$ to $\text{St}(m,n) $ so
\begin{equation}
	\text{grad}f(X) = P_X \text{grad}\bar{f}(X), \quad \text{where  } \text{ grad}\bar{f}(X) = 2AXN.
\end{equation}

The unbalanced orthogonal Procrustes problem consists of minimizing the function
\begin{equation}
	f:\text{St}(m,n) \rightarrow \mathbb{R}, \quad         X\mapsto f(X) = \| AX-B \|_F^2 ,
\end{equation}
on the Stiefel manifold $\text{St}(m,n) $, for given matrices $A \in \mathbb{R}^{l\times n }$ and $B \in \mathbb{R}^{l\times m }$ with $l \geq n$ and $l>m$, where $\| \cdot \|_F$ is the Frobenius norm $\| X \|_F^2 = \text{Trace}(X^\top X)=  \sum_{ij}{X_{ij}^2}$. If we define $\bar{f} : \mathbb{R}^{n\times m} \rightarrow \mathbb{R}$ via $X\mapsto \bar{f}(X) = \| AX-B \|_F^2,$ then $f$ is the restriction of $\bar{f}$ to $\text{St}(m,n) $ so
\begin{equation}
	 \text{grad}f(X) = P_X \text{grad}\bar{f}(X), \quad  \text{ where  } \text{ grad}\bar{f}(X) =  2A^\top (AX-B).
\end{equation}
Note that the special case where $n=m$ is the balanced orthogonal Procrustes problem. In this case, $\text{St}(m,n) = O(n)$ so $\|  AX  \|_F^2 = \|  A  \|_F^2$ and minimizing the function $f(X) = \| AX-B \|_F^2$ is replaced by the problem of maximizing $\text{Trace}(X^\top A^\top B)$ over $X\in O(n)$. A solution is then given by $X^* = UV^\top  $ where $B^\top A = U \Sigma V^\top $ is the Singular Value Decomposition of $B^\top A$ with square orthogonal matrices $U$ and $V$, and the solution is unique provided $B^\top A$ is nonsingular (see~\cite{Elden1999,Golub2013}).  \\

\subsection{Numerical Methods}

\subsubsection{\textbf{H}amiltonian \textbf{T}aylor \textbf{V}ariational \textbf{I}ntegrators (HTVIs)} 

HTVIs were first introduced in~\cite{ScShLe2017}. A discrete approximate Hamiltonian is constructed by approximating the flow map and the trajectory associated with the boundary values using a Taylor method, and approximating the integral by a quadrature rule. The Hamiltonian Taylor variational integrator is then generated by the discrete Hamilton's equations. More explicitly, Type II HTVIs are constructed as follows:

\begin{enumerate}[label=(\roman*)]
	\item Construct the $r$-order and $(r+1)$-order Taylor methods $\Psi_h^{(r)}$ and $\Psi_h^{(r+1)}$ approximating the exact time-$h$ flow map $\Phi _h : T^*Q \rightarrow T^*Q$.
	
	\item Approximate $p(0)=p_0$ by the solution $\tilde{p}_0$ of
	$ p_1 = \pi_{T^*Q} \circ \Psi_h^{(r)}(q_0,\tilde{p}_0) , $
	where $\pi_{T^*Q}:(q,p)\mapsto p$.
	
	\item Choose a quadrature rule of order $s$ with weights and nodes given by $(b_i,c_i)$ for $i=1,...,m$ and generate approximations $(q_{c_i},p_{c_i}) \approx (q(c_i h),p(c_i h))$ via $ (q_{c_i},p_{c_i})  = \Psi_{c_i h}^{(r)}(q_0,\tilde{p}_0).$
	
	\item Approximate $q_1$ via $ \tilde{q}_1 = \pi_{Q}  \circ \Psi_h^{(r+1)}(q_0,\tilde{p}_0),$ where $\pi_{Q}:(q,p)\mapsto q$.
	
	\item Use the continuous Legendre transform to obtain $\dot{q}_{c_i} = \frac{\partial H}{\partial p_{c_i}}$.
	
	\item Apply the quadrature rule to obtain the associated discrete right Hamiltonian
	
	$\label{HTVI_Hd}
		H_d^+(q_0,p_1) = p_1 \tilde{q}_1 - h \sum_{i=1}^{m}{b_i \left[  p_{c_i} \dot{q}_{c_i} - H(q_{c_i},p_{c_i})  \right]}.$
	
	\item The variational integrator is then defined by the discrete right Hamilton's equations.
	
\end{enumerate} 

\noindent Note that the following error analysis result concerning the order of accuracy of HTVIs was derived in~\cite{duruisseaux2020adaptive} (it can be extended to the constrained case via the strategy and results of Section~\ref{sec: Constrained Error Analysis}):
\begin{theorem} \label{HTVITheorem}
	If the Hamiltonian $H$ and its partial derivative $\frac{\partial H}{\partial p}$ are Lipschitz continuous in both variables, then $H_d^+(q_0,p_1)$ approximates $H_d^{+,E}(q_0,p_1)$ with at least order of accuracy $\min{(r+1,s)}$. 

	\noindent By Theorem 2.2 in~\cite{ScLe2017}, the associated discrete Hamiltonian map has the same order of accuracy.
\end{theorem}

\noindent In this paper, we will use the Direct approach and Adaptive approach $r=0$ Type II HTVIs constructed in~\cite{duruisseaux2020adaptive} based on the Direct and Adaptive discrete right Hamiltonians (respectively)
\begin{align}
	H_d^+(\bar{q}_0,\bar{r}_1;h) & = r_1^\top  q_0 + r_1^t q_0^t + h  \frac{p}{2(q_0^t)^{p+1}} r_1^\top  r_1   + hCp(q_0^t)^{2p-1} f(q_0)  + hr_1^t ,  \\
	H_d^+(\bar{q}_0,\bar{r}_1;h) & = r_1^\top  q_0 + r_1^t q_0^t + h \frac{p^2}{2\mathring{p} (q_0^t)^{p+\frac{\mathring{p}}{p}}  }   r_1^\top r_1  + hC \frac{p^2}{\mathring{p}} (q_0^t)^{2p-\frac{\mathring{p}}{p}} f(q_0) +   h \frac{p}{\mathring{p}}  (q_0^t)^{1-\frac{\mathring{p}}{p}}  r_1^t  .	
\end{align}

\begin{algorithm}[H] \label{Alg: HTVI}
	\DontPrintSemicolon
	
	\KwInput{A function  $f : \mathcal{Q} \rightarrow \mathbb{R}$, constants $C,h,p, \mathring{p} >0$,   $q_0^t , r_0^t \in \mathbb{R},$ and $(q_0 , r_0, \lambda_0)  \in T^*_{q_{0}} \mathcal{Q} \times \Lambda$.}
	
	 \textbf{while} convergence criterion is not met, \textbf{solve} the following system of equations:

\noln \begin{minipage}{0.455\textwidth}
	\hfill 
	
	\centering
	\textbf{Direct Approach} 
	
	\hfill 
	
	\hfill
	
	\centering
		$\begin{aligned}
		0 & =  r_{k+1} - r_k + hCp(q_k^t)^{2p-1} \nabla f(q_k) + \lambda_k ^\top   \nabla \mathcal{C}(q_k)\\
		0 & =  r_{k+1}^t - r_k^t - h \frac{p(p+1)}{2(q_k^t)^{p+2}} r_{k+1}^\top r_{k+1} \\ & \qquad  \qquad \quad + hCp(2p-1)(q_k^t)^{2p-2} f(q_k)\\
		0  &= q_{k+1} - q_k - h \frac{p}{(q_k^t)^{p+1}} r_{k+1}  \\
		0  & = q_{k+1}^t - q_k^t - h \\
		0 & = \mathcal{C}(q_{k+1})
	\end{aligned}$
\hfill 

\hfill 

\hfill

\end{minipage}
\hfill 
\vrule
\hfill 
\begin{minipage}{0.455\textwidth}
		\hfill

\centering
\textbf{Adaptive Approach}   

	\hfill 
	
\centering
	$\begin{aligned}
	0 & = r_{k+1} - r_k  + \frac{hCp^2}{\mathring{p}} (q_k^t)^{2p-\frac{\mathring{p}}{p}} \nabla f(q_k) +\lambda_k ^\top   \nabla \mathcal{C}(q_k) \\
	0 & = r_{k+1}^t - r_k^t  +  \frac{p \mathring{p} - 2p^3}{\mathring{p}}  h C (q_k^t)^{2p - \frac{\mathring{p}}{p} -1 } f(q_k) \\
	& 
	\qquad   + h\frac{p^3 +  p \mathring{p}}{2 \mathring{p} (q_k^t)^{p+\frac{\mathring{p}}{p} + 1}}  r_{k+1}^\top r_{k+1} + \frac{\mathring{p} - p}{ \mathring{p} (q_k^t)^{\frac{\mathring{p}}{p}}} h r_{k+1}^t  \\
	0 &= q_{k+1} - q_k - \frac{p^2}{\mathring{p} } h(q_k^t)^{-p-\frac{\mathring{p}}{p}} r_{k+1} \\
	0  & = q_{k+1}^t - q_k^t  - \frac{p}{\mathring{p} } h(q_k^t)^{1-\frac{\mathring{p}}{p}} \\
	0 & = \mathcal{C}(q_{k+1})
\end{aligned}$
\end{minipage}

%
	\caption{ Direct and Adaptive Hamiltonian Taylor variational integrators (HTVIs)}
\end{algorithm}

\subsubsection{Euler--Lagrange Simple Discretization} In~\cite{Duruisseaux2021Riemannian}, the $p$-Bregman Euler--Lagrange equations were rewritten as a first-order system of differential equations, for which a Riemannian version of a semi-implicit Euler scheme was applied to obtain the following algorithm:

\begin{algorithm}[H] \label{Alg: Semi-Implicit Euler}
	\DontPrintSemicolon
	
	\KwInput{A geodesically-convex ($\lambda = \zeta$) or $\alpha$-WQC ($\lambda = \zeta / \alpha$) function $f : \mathcal{Q} \rightarrow \mathbb{R}$. \\ \hspace{12.5mm} A retraction $\mathcal{R}$ from $T\mathcal{Q}$ to $\mathcal{Q}$, constants $C,h,p>0$, and $X_0 \in \mathcal{Q}$, $V_0 \in T_{X_0} \mathcal{Q}$. }
	
	\While{convergence criterion is not met}
	{
$b_k  \leftarrow 1 -  \frac{\lambda p+1}{k}, \quad c_k \leftarrow C p^2 (kh)^{p-2} $ \\
		\textbf{Version I}:	$a_k \leftarrow b_k V_k - hc_k \text{grad}f(X_k)$ \;
		\textbf{Version II}:	$a_k \leftarrow b_k V_k - h c_k \text{grad}f\left(\mathcal{R}_{X_k}(hb_k V_k)\right)$ \;
		
		$X_{k+1} \leftarrow \mathcal{R}_{X_k}(ha_k), \quad V_{k+1} \leftarrow \Gamma_{X_k}^{X_{k+1}} a_k$ \;
	} 
	\caption{Semi-Implicit Euler Integration of the $p$-Bregman Euler--Lagrange Equations }
\end{algorithm} 

Version I of Algorithm~\ref{Alg: Semi-Implicit Euler} corresponds to the usual update for the semi-implicit Euler scheme, while Version II is inspired by the reformulation of Nesterov's method from~\cite{Sutskever2013} that uses a corrected gradient $\nabla f(X_k +h b_kV_k)$ instead of the traditional gradient $\nabla f(X_k)$.

\subsubsection{\textbf{R}iemannian \textbf{G}radient \textbf{D}escent (RGD)} This is a generalization of Gradient Descent to the setting of Riemannian manifolds which involves the Riemannian gradient and a retraction.

\begin{algorithm}[H] \label{Alg: RGD}
	\DontPrintSemicolon

	\KwInput{A function $f : \mathcal{Q} \rightarrow \mathbb{R}$, a retraction $\mathcal{R}$ from $T\mathcal{Q}$ to $\mathcal{Q}$, $h>0$, and $X_0 \in  \mathcal{Q}$.}
	
	\While{convergence criterion is not met}
	{$X_{k+1} =  \mathcal{R}_{X_k} (-h \text{ grad} f(X_k))$} 
	\caption{Riemannian Gradient Descent (RGD)}
\end{algorithm}  

\hfill 

\subsection{Numerical Results}

It was noted in~\cite{Duruisseaux2021Riemannian} that although higher values of $p$ in Algorithm~\ref{Alg: Semi-Implicit Euler} result in provably faster rates of convergence, they also appear to be more prone to stability issues under numerical discretization, which can cause the numerical optimization algorithm to diverge. Numerical experiments in~\cite{duruisseaux2020adaptive} showed that in the normed vector space setting, geometric discretizations which respect the time-rescaling invariance and symplecticity of the Bregman Lagrangian and Hamiltonian flows were substantially less prone to these stability issues, and were therefore more robust, reliable, and computationally efficient. This was one of the motivations to develop time-adaptive Hamiltonian variational integrators for the Bregman Hamiltonians. Numerical experiments were conducted for the Rayleigh quotient minimization problem on $\mathbb{S}^{n-1}$, and for the generalized eigenvalue and Procrustes problems on the Stiefel manifold $\text{St}(m,n)$.

\begin{figure}[!hp]
	\centering
	\begin{minipage}[b]{0.9\textwidth}
		\includegraphics[width=\textwidth]{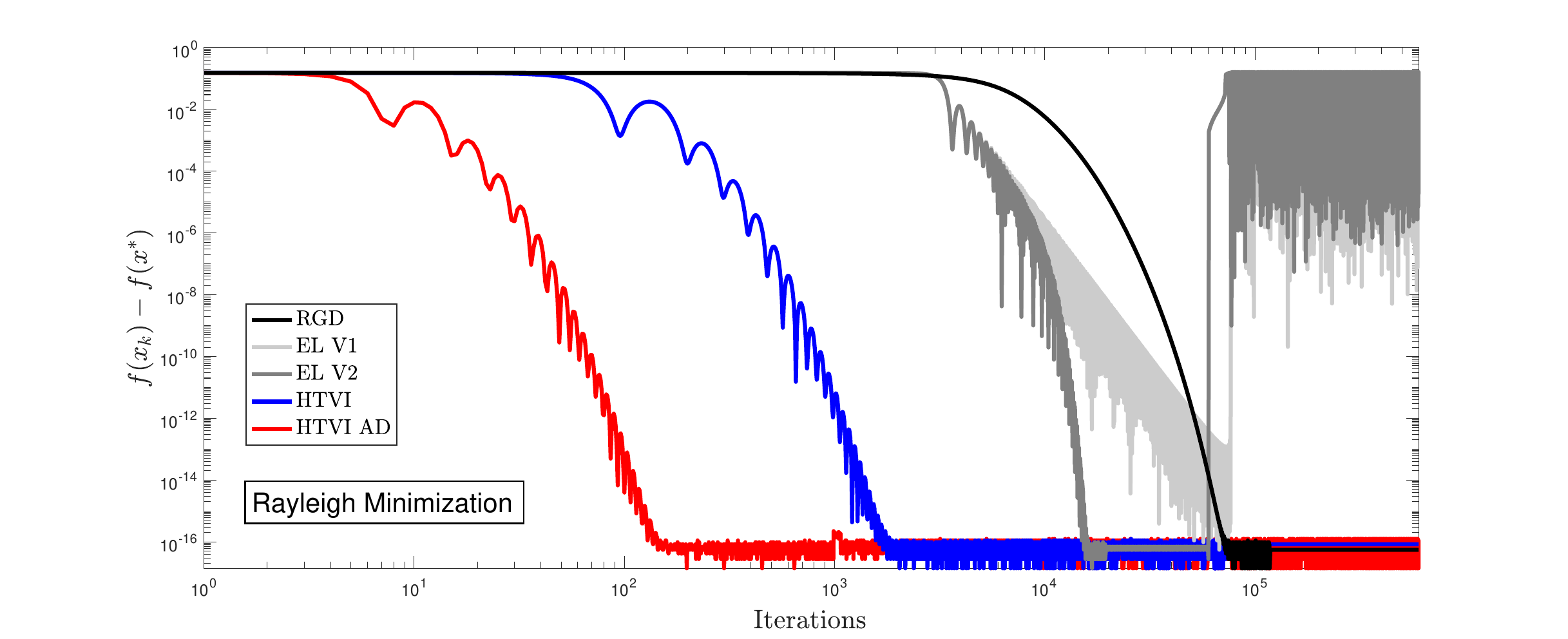}
	\end{minipage}
	\\
	\begin{minipage}[b]{0.9\textwidth}
		\includegraphics[width=\textwidth]{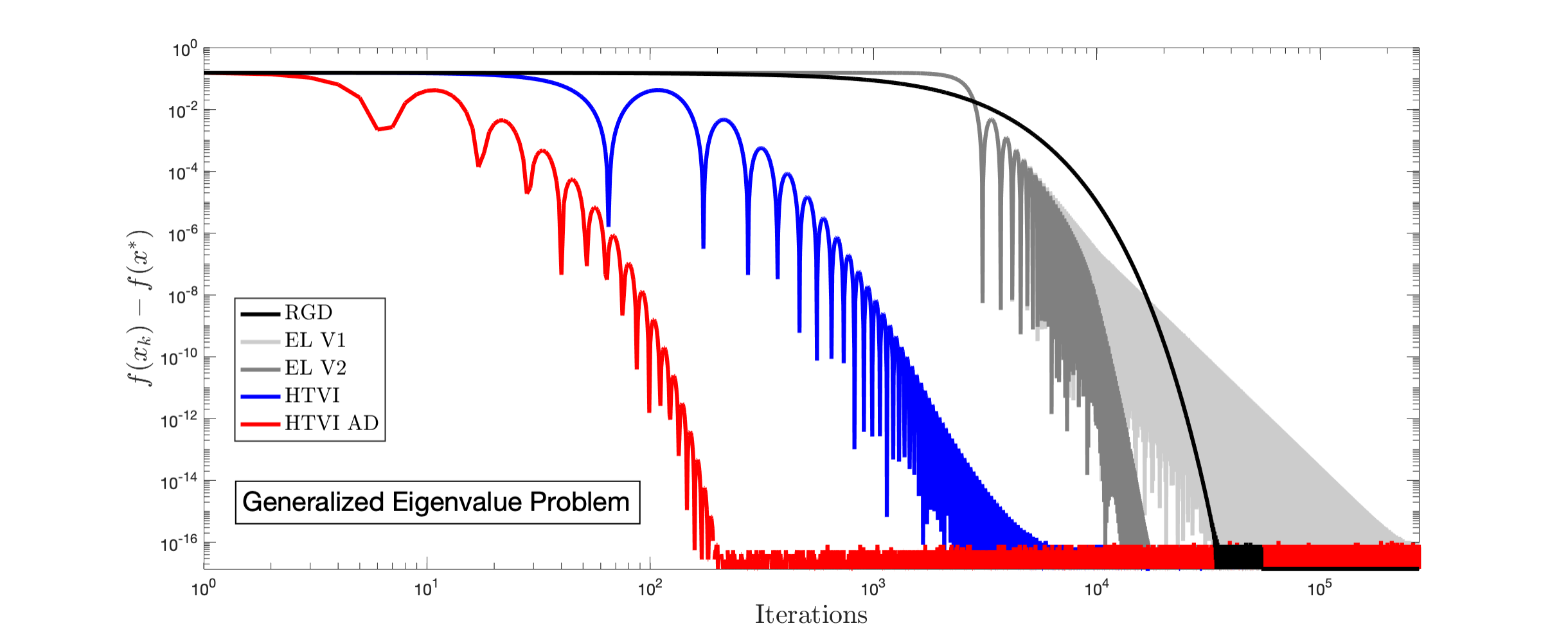}
	\end{minipage}
	\\
	\begin{minipage}[b]{0.9\textwidth}
		\includegraphics[width=\textwidth]{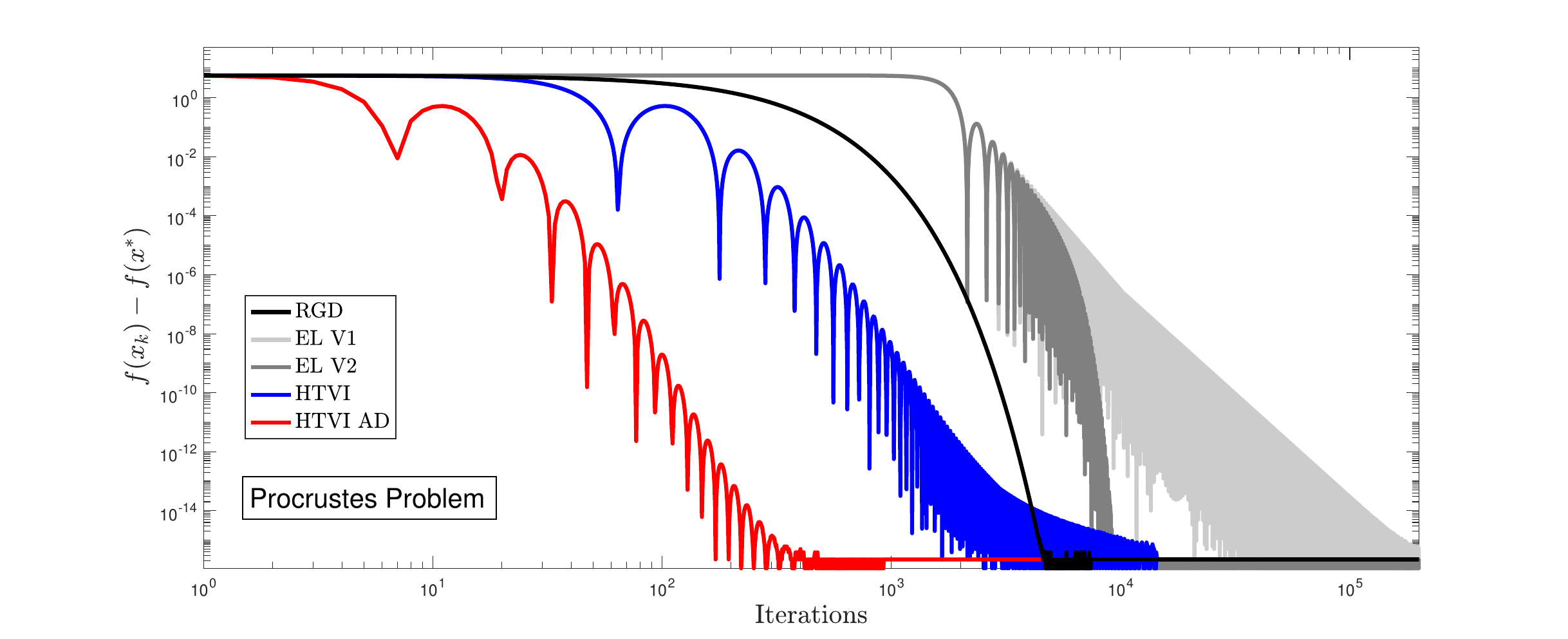}
	\end{minipage} 
	\caption
	{Comparison of the Direct and Adaptive (AD) Type II HTVIs with the Riemannian Gradient Descent (RGD) method and the Euler--Lagrange discretizations (EL V1 and EL V2) from~\cite{Duruisseaux2021Riemannian} with $p=6$ and the same timestep $h = 0.001$,  for the Rayleigh quotient minimization problem on the unit sphere $\mathbb{S}^{n-1}$, and for the generalized eigenvalue and Procrustes problems on the Stiefel manifold $\text{St}(m,n)$.  \label{fig: Stiefel} }
\end{figure}

The results from Figure~\ref{fig: Stiefel} show how the Hamiltonian Taylor variational integrators compare to the Euler--Lagrange discretizations from~\cite{Duruisseaux2021Riemannian} and the standard Riemannian gradient descent. Note that for certain instances of the Procrustes problem with certain initial values, all the algorithms converged to a local minimizer, and not the global minimizer, of the objective function. We can observe from Figure~\ref{fig: Stiefel} that for the same value of the timestep $h$, the Adaptive Hamiltonian variational integrator clearly outperforms its Direct counterpart, Riemannian gradient descent and the Euler--Lagrange discretizations in terms of number of iterations required. Furthermore, unlike the Euler--Lagrange discretizations (Algorithm~\ref{Alg: Semi-Implicit Euler}) and the Riemannian gradient descent (Algorithm~\ref{Alg: RGD}), the HTVI methods (Algorithm~\ref{Alg: HTVI}) do not require the use of retractions or parallel transports. Note that the Rayleigh minimization results indicate that the Euler--Lagrange discretizations suffer from stability issues leading to a loss of convergence, as the polynomially growing unbounded coefficient  $C p^2 (kh)^{p-2} $ is multiplied with $\text{grad} f$, so for this product to be bounded, the gradient has to decay to zero, but due to finite numerical precision, the gradient remains bounded away from zero, thereby causing the product to grow without bound. This issue can be resolved by adding a suitable upper bound to the coefficient $C p^2 (kh)^{p-2} $ in the updates, as can be seen both for the Euler--Lagrange discretizations and Hamiltonian variational integrators for the problems on $\text{St}(m,n)$.

However, the algorithms generated by these constrained Hamiltonian variational integrators are implicit, which can significantly increase the cost per iteration as the dimension of the problem becomes very large. In this case, it might be beneficial to consider other options using the unconstrained explicit Hamiltonian Taylor variational integrator, such as incorporating the constraints within the objective function as a penalty, although this might not constrain the solution trajectory to lie exactly on the manifold, or using projections if they can be computed efficiently and accurately for the Riemannian manifold of interest~\cite{duruisseaux2021projection}. Further, note that the implementation of the Hamiltonian variational integrators needs a very careful tuning of the various parameters at play, which may be challenging and thus also motivates the development of different methods. \\

\section{Conclusion}

Motivated by variational formulations of optimization problems on Riemannian manifolds, we first studied the relationship between the constrained Type I/II/III variational principles and the corresponding constrained Hamilton's or Euler--Lagrange equations both in continuous and discrete time, and derived variational error analysis results for the maps defined implicitly by the resulting discrete constrained equations. We then exploited these discrete constrained variational integrators and the variational formulation of accelerated optimization on Riemannian manifolds from~\cite{Duruisseaux2021Riemannian} to numerically solve the generalized eigenvalue and Procrustes problems on $\mathbb{S}^{n-1}$ and $\text{St}(m,n)$. 

The numerical experiments conducted in this paper corroborated the observation made for the vector space setting in~\cite{duruisseaux2020adaptive} that the Adaptive Hamiltonian variational integrator is significantly more efficient than the Direct Hamiltonian variational integrator, and that it can significantly outperform the Euler--Lagrange discretizations and Riemannian gradient descent, when its parameters are tuned carefully. Furthermore, it was noted that unlike the Euler--Lagrange discretizations and Riemannian gradient descent, the Hamiltonian algorithms did not require the use of retractions or parallel transports, which could be important when the problem considered lies on a Riemannian manifold for which it might not be possible to compute or approximate these objects efficiently. 

We noted however that tuning the parameters of these discrete constrained variational integrators can be challenging, and also that the resulting algorithms are implicit, which may significantly increase the cost per iteration as the dimension of the problem becomes very large, in which case it might be beneficial to consider using the unconstrained explicit HTVIs with projections~\cite{duruisseaux2021projection} or by incorporating the constraints within the objective function as a penalty. Moreover, although the Whitney and Nash Embedding Theorems~\cite{Whitney1944_2,Whitney1944_1,Nash1956} imply that there is no loss of generality when studying Riemannian manifolds only as submanifolds of Euclidean spaces, there are limitations to the constrained integration strategy based on embeddings presented in this paper, and an approach intrinsically defined on Riemannian manifolds would be desirable. Indeed, the embedding approach usually leads to higher-dimensional computations, and requires an effective way of constructing the embedding or a natural way of writing down equations that constrain the problem and the numerical solutions to the Riemannian manifold. Furthermore, most results in Riemannian geometry or results concerning specific Riemannian manifolds are proven from an intrinsic perspective because the embedding approach tends to flood intrinsic geometric properties of the manifold with superfluous information coming from the additional dimensions of the Euclidean space. This motivates the development of intrinsic methods that would exploit the symmetries and geometric properties of the manifold and of the problem at hand.

Developing an intrinsic extension of Hamiltonian variational integrators to manifolds will require some additional work, since the current approach involves Type II/III generating functions $H_d^+(q_k, p_{k+1})$, $H_d^-(p_k, q_{k+1})$, which depend on the position at one boundary point, and the momentum at the other boundary point. However, this does not make intrinsic sense on a manifold, since one needs the base point in order to specify the corresponding cotangent space, and one should ideally consider a Hamiltonian variational integrator construction based on discrete Dirac mechanics~\cite{LeOh2008}, which would yield a generating function $E_d^+(q_k, q_{k+1}, p_{k+1})$, $E_d^-(q_k, p_k, q_{k+1})$, that depends on the position at both boundary points and the momentum at one of the boundary points. This approach can be viewed as a discretization of the generalized energy $E(q,v,p)=\langle p,v\rangle - L(q,v)$, in contrast to the Hamiltonian $H(q,p)=\ext_{v}\langle p,v\rangle - L(q,v)=\left.\langle p,v\rangle - L(q,v)\right|_{p=\frac{\partial L}{\partial v}}$. On the other hand, the formulation of Lagrangian variational integrators presented in the introduction of Section~\ref{section: Constrained Variational Lagrangian} makes sense intrinsically on manifolds, and a framework for variable time-stepping in Lagrangian variational integration was introduced in our subsequent work~\cite{Duruisseaux2022Lagrangian} to design intrinsic time-adaptive Lagrangian variational integrators for accelerated optimization on Riemannian manifolds.

It would also be interesting to extend the proposed approach to the problem of optimization with nonintegrable constraints, which naturally leads to the question of whether vakonomic or nonholonomic mechanics is the appropriate description~\cite{Cortes2002}. In the context of optimization with nonintegrable constraints, the relevant extension will likely involve vakonomic variational integrators~\cite{Benito2005,Jimenez2012}. However, it would be interesting to relate the methods introduced in this paper to the existing work on variational integrators applied to optimal control problems~\cite{Junge2005,Leon2007}, and the discrete optimal control of nonholonomic dynamical systems would likely require a combination of the methods described here and nonholonomic integrators~\cite{Cortes2001,Leon2003, Fedorov2005, McLachlan2006}.


\section*{Acknowledgments} 

The authors were supported in part by NSF under grants DMS-1411792, DMS-1345013, DMS-1813635, CCF-2112665, by AFOSR under grant FA9550-18-1-0288, and by the DoD under grant 13106725 (Newton Award for Transformative Ideas during the COVID-19 Pandemic).

\section*{Data Availability Statement}
The datasets generated during and/or analyzed during the current study are available from the corresponding author on reasonable request.

  \appendix

\section{Proofs of Theorems for Constrained Variational Mechanics}

\subsection{Proof of Theorem~\ref{theorem: Lagrangian Action and Equations}} \label{Appendix: Type I Action}

\begin{theorem} 
	Consider the constrained action functional $\mathfrak{S} : C^2([0,T],\mathcal{Q} \times \Lambda) \rightarrow \mathbb{R}$ given by
	\begin{equation}
		\mathfrak{S} (q(\cdot), \lambda(\cdot)) = \int_{0}^{T}{ \left[ L(q(t),\dot{q}(t)) - \langle \lambda(t) ,  \mathcal{C}(q(t)) \rangle  \right] dt}.
	\end{equation}
	The condition that $\mathfrak{S} (q(\cdot), \lambda(\cdot)) $ is stationary with respect to the boundary conditions $\delta q(0) = 0$ and $\delta q(T) =0$ is equivalent to $(q(\cdot),\lambda(\cdot))$ satisfying the constrained Euler--Lagrange equations
	\begin{equation} \label{ELConstEqsAppendix}
		\frac{\partial L}{\partial q} - \frac{d}{dt} \frac{\partial L}{\partial \dot{q}}  =  \langle \lambda , \nabla \mathcal{C} (q) \rangle, \qquad \mathcal{C}(q) =0 . 
	\end{equation}
\small
	\proof{Computing the variation of $\mathfrak{S}$ yields
		\begin{align*} 
			\delta \mathfrak{S} & = \int_{0}^{T}{ \left[ \frac{\partial L}{\partial q}(q(t),\dot{q}(t)) \delta q(t) + \frac{\partial L}{\partial \dot{q}}(q(t),\dot{q}(t)) \delta \dot{q}(t)\right] dt}   -  \int_{0}^{T}{ \left[  \langle \lambda (t) , \nabla \mathcal{C} (q(t)) \delta q(t) \rangle  + \langle \delta \lambda (t) , \mathcal{C}(q(t)) \rangle  \right] dt } .
		\end{align*}  
		Using integration by parts and the boundary conditions $\delta q(0) = 0$ and $\delta q (T) = 0,$ we get
		\begin{align*}
			\delta \mathfrak{S} & =\int_{0}^{T}{ \left[ \frac{\partial L}{\partial q}(q(t),\dot{q}(t)) - \frac{d}{dt} \frac{\partial L}{\partial \dot{q}}(q(t),\dot{q}(t))  -  \langle \lambda (t) , \nabla \mathcal{C} (q(t))  \rangle  \right] \delta q(t) dt } - \int_{0}^{T}{  \langle \delta \lambda (t) , \mathcal{C}(q(t)) \rangle  dt }  .
		\end{align*}  
		Now, if $\delta \mathfrak{S} = 0$, then the fundamental theorem of the calculus of variations~\cite{Ar1989} yields the constrained Euler--Lagrange equations \eqref{ELConstEqsAppendix}. Conversely, if $(q,\lambda )$ satisfies the constrained Euler--Lagrange equations \eqref{ELConstEqsAppendix}, then the integrand vanishes and $\delta \mathfrak{S} = 0$. \qed} \\
\end{theorem}

\subsection{Proof of Theorem~\ref{theorem: Action and Equations II}} \label{Appendix: Type II Action}

\begin{theorem} 
	Consider the constrained action functional $\mathfrak{S} : C^2([0,T],T^*\mathcal{Q} \times \Lambda) \rightarrow \mathbb{R}$ given by
	\begin{equation}
		\mathfrak{S} (q(\cdot),p(\cdot), \lambda(\cdot)) = p(T)q(T) - \int_{0}^{T}{ \left[ p(t) \dot{q}(t) - H(q(t),p(t)) - \langle \lambda(t) ,  \mathcal{C}(q(t)) \rangle  \right] dt}.
	\end{equation}
	The condition that $\mathfrak{S} (q(\cdot),p(\cdot), \lambda(\cdot)) $ is stationary with respect to the boundary conditions $\delta q(0) = 0$ and $\delta p(T) =0$ is equivalent to $(q(\cdot),p(\cdot),\lambda(\cdot))$ satisfying Hamilton's canonical constrained equations
	\begin{equation} \label{HamiltonConstEqsAppendix}
		\dot{q} = \frac{\partial H}{\partial p} (q,p) , \qquad \dot{p} = -\frac{\partial H}{\partial q} (q,p) - \langle \lambda,  \nabla \mathcal{C} (q) \rangle, \qquad \mathcal{C}(q) = 0 . 
	\end{equation}
\small 
	\proof{Computing the variation of $\mathfrak{S}$ yields
		\begin{align} \label{ProofActionEq1}
			\delta \mathfrak{S} & = q(T) \delta p(T) + p(T) \delta q (T)  + \int_{0}^{T}{ \left[  \langle \lambda (t) , \nabla \mathcal{C} (q(t)) \delta q(t) \rangle  + \langle \delta \lambda (t) , \mathcal{C}(q(t)) \rangle  \right] dt }   \nonumber \\ & \qquad \qquad  \quad - \int_{0}^{T}{ \left[ \dot{q} (t) \delta p(t) + p(t) \delta \dot{q}(t) - \frac{\partial H}{\partial q}(q(t),p(t)) \delta q(t) - \frac{\partial H}{\partial p}(q(t),p(t)) \delta p(t) \right] dt }. 
		\end{align}  
		Using integration by parts and the boundary conditions $\delta q(0) = 0$ and $\delta p (T) = 0,$ we get
		\begin{align*}
			\delta \mathfrak{S} & = q(T) \delta p(T) + p(T) \delta q(T) - p(T) \delta q(T)  + p(0) \delta q(0) +\int_{0}^{T}{  \langle \delta \lambda (t) , \mathcal{C}(q(t)) \rangle  dt }   \\ 
			&  \qquad \quad+ \int_{0}^{T}{ \left[ \dot{p}(t) + \frac{\partial H}{\partial q}(q(t),p(t))  + \langle \lambda(t) , \nabla \mathcal{C} (q(t))  \rangle \right] \delta q(t) dt } + \int_{0}^{T}{ \left[ \frac{\partial H}{\partial p}(q(t),p(t))  -\dot{q} (t) \right]  \delta p(t) dt }  \\
			& =  \int_{0}^{T}{ \left[ \dot{p}(t) + \frac{\partial H}{\partial q}(q(t),p(t))  +  \langle \lambda(t) , \nabla \mathcal{C} (q(t)) \rangle \right] \delta q(t) dt } + \int_{0}^{T}{ \left[ \frac{\partial H}{\partial p}(q(t),p(t))  -\dot{q} (t) \right]  \delta p(t) dt }   \\ & \quad \qquad + \int_{0}^{T}{  \langle \delta \lambda (t) , \mathcal{C}(q(t)) \rangle  dt }   .
		\end{align*}  
		Now, if $\delta \mathfrak{S} = 0$, then the fundamental theorem of the calculus of variations~\cite{Ar1989} yields Hamilton's constrained equations \eqref{HamiltonConstEqsAppendix}. Conversely, if $(q,p,\lambda )$ satisfies Hamilton's constrained equations \eqref{HamiltonConstEqsAppendix}, then the integrand vanishes and $\delta \mathfrak{S} = 0$. \qed} \\
\end{theorem}

\normalsize
\subsection{Proof of Theorem~\ref{theorem: Action and Equations III}} \label{Appendix: Type III Action}

\begin{theorem} 
	Consider the constrained action functional $\mathfrak{S} : C^2([0,T],T^*\mathcal{Q} \times \Lambda) \rightarrow \mathbb{R}$ given by
	\begin{equation}
		\mathfrak{S} (q(\cdot),p(\cdot), \lambda(\cdot)) = - p(0)q(0) - \int_{0}^{T}{ \left[ p(t) \dot{q}(t) - H(q(t),p(t)) - \langle \lambda(t) ,  \mathcal{C}(q(t)) \rangle \right] dt}.
	\end{equation}
	The condition that $\mathfrak{S} (q(\cdot),p(\cdot), \lambda(\cdot)) $ is stationary with respect to the boundary conditions $\delta q(T) = 0$ and $\delta p(0) =0$ is equivalent to $(q(\cdot),p(\cdot),\lambda(\cdot))$ satisfying Hamilton's canonical constrained equations
	\begin{equation} \label{HamiltonConstEqsIII Appendix}
		\dot{q} = \frac{\partial H}{\partial p} (q,p) , \qquad \dot{p} = -\frac{\partial H}{\partial q} (q,p) - \langle \lambda,  \nabla \mathcal{C} (q) \rangle,  \qquad \mathcal{C}(q) = 0 .
	\end{equation}
\small  \proof{ The proof is almost identical to that of Theorem~\ref{theorem: Action and Equations II}. We compute the variation of $\mathfrak{S}$ as before and get equation \eqref{ProofActionEq1} except that the term $\left( q(T) \delta p(T) + p(T) \delta q (T) \right)$ is replaced by $\left(- q(0) \delta p(0) - p(0) \delta q (0)\right)$. As before, integration by parts and the boundary conditions $\delta q(T) = 0$ and $\delta p(0) =0$ yield
		\begin{align*}
			\delta \mathfrak{S} & =  \int_{0}^{T}{ \left[ \dot{p}(t) + \frac{\partial H}{\partial q}(q(t),p(t))  + \langle \lambda(t) , \nabla \mathcal{C} (q(t)) \rangle \right] \delta q(t) dt } + \int_{0}^{T}{ \left[ \frac{\partial H}{\partial p}(q(t),p(t))  -\dot{q} (t) \right]  \delta p(t) dt }  \\ & \quad \qquad + \int_{0}^{T}{  \langle \delta \lambda (t) , \mathcal{C}(q(t)) \rangle     dt }   .
		\end{align*}  
		Then, if $\delta \mathfrak{S} = 0$, then the fundamental theorem of the calculus of variations~\cite{Ar1989} yields Hamilton's constrained equations \eqref{HamiltonConstEqsIII Appendix}. Conversely, if $(q,p,\lambda )$ satisfies Hamilton's constrained equations \eqref{HamiltonConstEqsIII Appendix}, then the integrand vanishes and $\delta \mathfrak{S} = 0$. \qed} \\
\end{theorem}

\subsection{Proof of Theorem~\ref{theorem: Flow Map and Equations I}} \label{Appendix: Type I Flow}

\begin{theorem} 
	The exact time-$T$ flow map of Hamilton's equations $(q_0,p_0) \mapsto (q_T,p_T)$ is implicitly given by the following relations:
	\begin{equation}
		D_1 \mathcal{S} (q_0,q_T) = -  \frac{\partial L}{ \partial \dot{q}}(q_0,\dot{q}(0)), \qquad  D_2 \mathcal{S} (q_0,q_T) =  \frac{\partial L}{ \partial \dot{q}}(q_T,\dot{q}(T))  .
	\end{equation}
	Thus, $ \mathcal{S} (q_0,q_T)$ is a Type I generating function that generates the exact flow of the constrained Euler--Lagrange equations \eqref{ELConstEqs}. 
	\small 	\proof{Using integration by parts and simplifying gives
		\begin{align*}
			\frac{ \partial \mathcal{S}}{ \partial q_0 }(q_0,q_T) & =  \int_{0}^{T}{ \left[ \frac{\partial q(t)}{\partial q_0}  \frac{\partial L}{\partial q}(q(t),\dot{q}(t)) + \frac{\partial \dot{q}(t)}{\partial q_0}  \frac{\partial L}{\partial \dot{q}}(q(t),\dot{q}(t)) \right] dt }   \\ & \qquad  \qquad \qquad   -  \int_{0}^{T}{ \left[ \langle \lambda(t) , \frac{\partial q(t)}{\partial q_0} \nabla \mathcal{C} (q(t)) \rangle +  \langle \frac{\partial \lambda (t) }{\partial q_0} , \mathcal{C}(q(t)) \rangle \right] dt}  \\
			& =  \int_{0}^{T}{   \frac{\partial q(t)}{\partial q_0}  \left( \frac{\partial L}{\partial q}(q(t),\dot{q}(t)) - \frac{d}{dt} \frac{\partial L}{\partial \dot{q}}(q(t),\dot{q}(t))  -  \langle \lambda (t) , \nabla \mathcal{C} (q(t))  \rangle    \right)  dt}  \\ &  \qquad \qquad   \qquad -   \int_{0}^{T}{  \langle \frac{\partial \lambda (t) }{\partial q_0} , \mathcal{C}(q(t)) \rangle dt} -  \frac{\partial L}{ \partial \dot{q}}(q(0),\dot{q}(0)) ,
		\end{align*}
		\begin{align*}
			\frac{ \partial \mathcal{S}}{ \partial q_T }(q_0,q_T) & =  \int_{0}^{T}{ \left[ \frac{\partial q(t)}{\partial q_T}  \frac{\partial L}{\partial q}(q(t),\dot{q}(t)) + \frac{\partial \dot{q}(t)}{\partial q_T}  \frac{\partial L}{\partial \dot{q}}(q(t),\dot{q}(t)) \right] dt }  \\ & \qquad  \qquad \qquad    -  \int_{0}^{T}{ \left[ \langle \lambda(t) , \frac{\partial q(t)}{\partial q_T} \nabla \mathcal{C} (q(t)) \rangle +  \langle \frac{\partial \lambda (t) }{\partial q_T} , \mathcal{C}(q(t)) \rangle \right] dt}  \\
			& =  \int_{0}^{T}{   \frac{\partial q(t)}{\partial q_T}  \left( \frac{\partial L}{\partial q}(q(t),\dot{q}(t)) - \frac{d}{dt} \frac{\partial L}{\partial \dot{q}}(q(t),\dot{q}(t))  -  \langle \lambda (t) , \nabla \mathcal{C} (q(t))  \rangle    \right)  dt}  \\ &  \qquad  \qquad  \qquad    -   \int_{0}^{T}{  \langle \frac{\partial \lambda (t) }{\partial q_T} , \mathcal{C}(q(t)) \rangle dt} + \frac{\partial L}{ \partial \dot{q}}(q(T),\dot{q}(T)).
		\end{align*}
		By Theorem~\ref{theorem: Lagrangian Action and Equations}, the extremum of the action functional $\mathfrak{S}$ is achieved when $(q,\lambda)$ satisfies the constrained Euler--Lagrange equations \eqref{ELConstEqs}, so we get $  D_1 \mathcal{S} (q_0,q_T) = -  \frac{\partial L}{ \partial \dot{q}}(q_0,\dot{q}(0))  $ and $  D_2 \mathcal{S} (q_0,q_T) =  \frac{\partial L}{ \partial \dot{q}}(q_T,\dot{q}(T))  $.  \qed \\
	}
\end{theorem}

\subsection{Proof of Theorem~\ref{theorem: Flow Map and Equations II}} \label{Appendix: Type II Flow}

\begin{theorem} 
	The exact time-$T$ flow map of Hamilton's equations $(q_0,p_0) \mapsto (q_T,p_T)$ is implicitly given by the following relations:
	\begin{equation}
		q_T = D_2 \mathcal{S} (q_0,p_T), \qquad p_0 = D_1 \mathcal{S} (q_0,p_T).
	\end{equation}
	In particular, $ \mathcal{S} (q_0,p_T)$ is a Type II generating function that generates the exact flow of Hamilton's constrained equations \eqref{HamiltonConstEqs}. 
	
\small 	\proof{Using integration by parts and simplifying gives
		\begin{align*}
			\frac{ \partial \mathcal{S}}{ \partial q_0 }(q_0,p_T) & = \frac{\partial q_T}{\partial q_0} p_T  + \int_{0}^{T}{ \left[ \langle \lambda(t) , \frac{\partial q(t)}{\partial q_0} \nabla \mathcal{C} (q(t)) \rangle +  \langle \frac{\partial \lambda (t) }{\partial q_0} , \mathcal{C}(q(t)) \rangle \right] dt} \\
			& \qquad \qquad - \int_{0}^{T}{ \left[ \frac{\partial p(t)}{\partial q_0} \dot{q}(t) + \frac{\partial \dot{q}(t)}{\partial q_0} p(t) - \frac{\partial q(t)}{\partial q_0} \frac{\partial H}{\partial q}(q(t),p(t))  - \frac{\partial p(t)}{\partial q_0} \frac{\partial H}{\partial p}(q(t),p(t))  \right] dt }   \\
			& = p_0 + \int_{0}^{T}{   \frac{\partial q(t)}{\partial q_0}  \left( \dot{p}(t) +   \frac{\partial H}{\partial q}(q(t),p(t)) + \langle \lambda(t) ,  \nabla \mathcal{C} (q(t)) \rangle   \right)  dt} 
			\\
			& \qquad \qquad  -  \int_{0}^{T}{  \frac{\partial p(t)}{\partial q_0} \left( \dot{q}(t) - \frac{\partial H}{\partial p}(q(t),p(t)) \right)  dt}  + \int_{0}^{T}{  \langle \frac{\partial \lambda (t) }{\partial q_0} , \mathcal{C}(q(t)) \rangle dt} ,
		\end{align*}
		\begin{align*}
			\frac{ \partial \mathcal{S}}{ \partial p_T }(q_0,p_T) & = q_T + \frac{\partial q_T}{\partial p_T} p_T + \int_{0}^{T}{ \left[ \langle \lambda(t) , \frac{\partial q(t)}{\partial p_T} \nabla \mathcal{C} (q(t)) \rangle +  \langle \frac{\partial \lambda (t) }{\partial p_T} , \mathcal{C}(q(t)) \rangle \right] dt}  \\
			& \qquad  \qquad - \int_{0}^{T}{ \left[ \frac{\partial p(t)}{\partial p_T} \dot{q}(t) + \frac{\partial \dot{q}(t)}{\partial p_T} p(t) - \frac{\partial q(t)}{\partial p_T} \frac{\partial H}{\partial q}(q(t),p(t))  - \frac{\partial p(t)}{\partial p_T} \frac{\partial H}{\partial p}(q(t),p(t))  \right] dt }   \\
			& = q_T + \int_{0}^{T}{   \frac{\partial q(t)}{\partial p_T}  \left( \dot{p}(t) +   \frac{\partial H}{\partial q}(q(t),p(t)) + \langle \lambda(t) , \nabla \mathcal{C} (q(t)) \rangle   \right)  dt} 
			\\
			& \qquad  \qquad -  \int_{0}^{T}{  \frac{\partial p(t)}{\partial p_T} \left( \dot{q}(t) - \frac{\partial H}{\partial p}(q(t),p(t)) \right)  dt}  + \int_{0}^{T}{  \langle \frac{\partial \lambda (t) }{\partial p_T} , \mathcal{C}(q(t)) \rangle  dt} .
		\end{align*}
		By Theorem~\ref{theorem: Action and Equations II}, the extremum of the action functional $\mathfrak{S}$ is achieved when the curve $(q,p,\lambda)$ satisfies Hamilton's constrained equations \eqref{HamiltonConstEqs}, so the integrands vanish, and thus  $p_0 = \frac{\partial \mathcal{S} }{ \partial q_0} (q_0,p_T) = D_1 \mathcal{S} (q_0,p_T)$ and $q_T = \frac{\partial \mathcal{S} }{ \partial p_T} (q_0,p_T) = D_2 \mathcal{S} (q_0,p_T)$. \qed \\
	}
\end{theorem}

\subsection{Proof of Theorem~\ref{theorem: Flow Map and Equations III}} \label{Appendix: Type III Flow}

\begin{theorem} 
	The exact time-$T$ flow map of Hamilton's equations $(q_0,p_0) \mapsto (q_T,p_T)$ is implicitly given by the following relations:
	\begin{equation}
		q_0 = -  D_2 \mathcal{S} (q_T,p_0), \qquad p_T = - D_1 \mathcal{S} (q_T,p_0).
	\end{equation}
	In particular, $ \mathcal{S} (q_T,p_0)$ is a Type III generating function that generates the exact flow of Hamilton's constrained equations \eqref{HamiltonConstEqsIII}. 
	
	\small \proof{Integrating by parts and simplifying yields
		\begin{align*}
			\frac{ \partial \mathcal{S}}{ \partial q_T }(q_T,p_0) & = - \frac{\partial q_0}{\partial q_T} p_0  + \int_{0}^{T}{ \left[ \langle \lambda(t) , \frac{\partial q(t)}{\partial q_T} \nabla \mathcal{C} (q(t)) \rangle +  \langle \frac{\partial \lambda (t) }{\partial q_T} , \mathcal{C}(q(t)) \rangle \right] dt} \\
			& \qquad  \qquad- \int_{0}^{T}{ \left[ \frac{\partial p(t)}{\partial q_T} \dot{q}(t) + \frac{\partial \dot{q}(t)}{\partial q_T} p(t) - \frac{\partial q(t)}{\partial q_T} \frac{\partial H}{\partial q}(q(t),p(t))  - \frac{\partial p(t)}{\partial q_T} \frac{\partial H}{\partial p}(q(t),p(t))  \right] dt }   \\
			& = - p_T + \int_{0}^{T}{   \frac{\partial q(t)}{\partial q_T}  \left( \dot{p}(t) +   \frac{\partial H}{\partial q}(q(t),p(t)) + \langle \lambda(t) ,  \nabla \mathcal{C} (q(t)) \rangle    \right)  dt} 
			\\
			& \qquad  \qquad -  \int_{0}^{T}{  \frac{\partial p(t)}{\partial q_T} \left( \dot{q}(t) - \frac{\partial H}{\partial p}(q(t),p(t)) \right)  dt}  + \int_{0}^{T}{  \langle \frac{\partial \lambda (t) }{\partial q_T} , \mathcal{C}(q(t)) \rangle dt} ,
		\end{align*} 
		\begin{align*}
			\frac{ \partial \mathcal{S}}{ \partial p_0 }(q_T,p_0) & = - q_0 -  \frac{\partial q_0}{\partial p_0} p_0 +\int_{0}^{T}{ \left[ \langle \lambda(t) , \frac{\partial q(t)}{\partial p_0} \nabla \mathcal{C} (q(t)) \rangle +  \langle \frac{\partial \lambda (t) }{\partial p_0} , \mathcal{C}(q(t)) \rangle \right] dt}\\
			& \qquad \qquad - \int_{0}^{T}{ \left[ \frac{\partial p(t)}{\partial p_0} \dot{q}(t) + \frac{\partial \dot{q}(t)}{\partial p_0} p(t) - \frac{\partial q(t)}{\partial p_0} \frac{\partial H}{\partial q}(q(t),p(t))  - \frac{\partial p(t)}{\partial p_0} \frac{\partial H}{\partial p}(q(t),p(t))  \right] dt }   \\
			& = -q_0 + \int_{0}^{T}{   \frac{\partial q(t)}{\partial p_0}  \left( \dot{p}(t) +   \frac{\partial H}{\partial q}(q(t),p(t)) + \langle \lambda(t) ,  \nabla \mathcal{C} (q(t)) \rangle    \right)  dt} 
			\\
			& \qquad \qquad -  \int_{0}^{T}{  \frac{\partial p(t)}{\partial p_0} \left( \dot{q}(t) - \frac{\partial H}{\partial p}(q(t),p(t)) \right)  dt}  + \int_{0}^{T}{  \langle \frac{\partial \lambda (t) }{\partial p_0} , \mathcal{C}(q(t)) \rangle dt} .
		\end{align*}
		By Theorem~\ref{theorem: Action and Equations III}, the extremum of the action $\mathfrak{S}$ is achieved when the curve $(q,p,\lambda)$ satisfies Hamilton's constrained equations \eqref{HamiltonConstEqsIII}, so the integrands vanish, and thus $p_T = - \frac{\partial \mathcal{S} }{ \partial q_T} (q_T,p_0) = - D_1 \mathcal{S} (q_T,p_0)$ and $q_0 = - \frac{\partial \mathcal{S} }{ \partial p_0} (q_T,p_0) = - D_2 \mathcal{S} (q_T,p_0)$. \qed \\
	}
\end{theorem}

\subsection{Proof of Theorem~\ref{theorem: Discrete Variational Principle I }} \label{Appendix: Type I Discrete Variational}

\begin{theorem}  The Type I discrete Hamilton's variational principles 
	\begin{equation}
		\delta \mathfrak{S}_d^{\pm} \left(\{  (q_k, \lambda _k) \}_{k=0}^{N} \right)  = 0
	\end{equation}
	are equivalent to the discrete constrained Euler--Lagrange equations
	\begin{equation} \label{Discrete Constrained EL equationsAppendix}
		D_1 L_d(q_k,q_{k+1})  + D_2 L_d(q_{k-1},q_{k}) = \langle \lambda_k ,\nabla \mathcal{C}(q_k) \rangle, \quad  \mathcal{C}(q_{k})  =0,
	\end{equation} 
	where $L_d(q_{k},q_{k+1}) $ is defined via equation \eqref{Discrete Lagrangian}. 
\small	\proof{Using the fact that $\delta q_0 =0$ and $\delta q_N =0$, we have
		\begin{align*}
			\delta \mathfrak{S}_d^-  & = \delta \left(   \sum_{k=0}^{N-1}{\left[  L_d(q_k,q_{k+1}) - \langle \lambda_k, \mathcal{C}(q_k)  \rangle \right]} \right) \\ & =  \sum_{k=0}^{N-1}{\left[  D_1 L_d(q_k,q_{k+1}) \delta q_k + D_2 L_d(q_k,q_{k+1}) \delta q_{k+1}  \right]}  - \sum_{k=0}^{N-1}{\left( \langle \lambda_k ,\nabla \mathcal{C}(q_k) \delta q_k \rangle    + \langle \delta \lambda_k , \mathcal{C}(q_k) \rangle  \right)} \\ 
			& =  \sum_{k=1}^{N-1}{\left[  D_1 L_d(q_k,q_{k+1})  + D_2 L_d(q_{k-1},q_{k}) - \langle \lambda_k ,\nabla \mathcal{C}(q_k) \rangle \right] \delta q_{k} }  - \sum_{k=0}^{N-1}{\langle \delta \lambda_k , \mathcal{C}(q_k) \rangle}, \\
		\delta \mathfrak{S}_d^+  & = \delta \left(   \sum_{k=0}^{N-1}{\left[  L_d(q_k,q_{k+1}) - \langle \lambda_{k+1}, \mathcal{C}(q_{k+1})  \rangle \right]} \right) \\ & = \sum_{k=0}^{N-1}{\left[  D_1 L_d(q_k,q_{k+1}) \delta q_k + D_2 L_d(q_k,q_{k+1}) \delta q_{k+1}  \right]}  - \sum_{k=0}^{N-1}{\left( \langle \lambda_{k+1} ,\nabla \mathcal{C}(q_{k+1}) \delta q_{k+1} \rangle    + \langle \delta \lambda_{k+1} , \mathcal{C}(q_{k+1}) \rangle  \right)} \\ 
		& =  \sum_{k=1}^{N-1}{\left[  D_1 L_d(q_k,q_{k+1})  + D_2 L_d(q_{k-1},q_{k}) - \langle \lambda_k ,\nabla \mathcal{C}(q_k) \rangle \right] \delta q_{k} }  - \sum_{k=0}^{N-1}{\langle \delta \lambda_{k+1} , \mathcal{C}(q_{k+1}) \rangle}. 
	\end{align*}
		If the discrete constrained Euler--Lagrange equations \eqref{Discrete Constrained EL equationsAppendix} are satisfied, then each term vanishes and $\delta \mathfrak{S}_d^{\pm} = 0$. Conversely, if $\delta \mathfrak{S}_d^{\pm} = 0$, then a discrete fundamental theorem of the calculus of variations yields the discrete constrained Euler--Lagrange equations  \eqref{Discrete Constrained EL equationsAppendix}. \qed \\
	}
\end{theorem}

\subsection{Proof of Theorem~\ref{theorem: Discrete Variational Principle II}} \label{Appendix: Type II Discrete Variational}

	\begin{theorem} 
	The Type II discrete Hamilton's phase space variational principle
	\begin{equation}
		\delta \mathfrak{S}_d^{+} \left(\{  (q_k,p_k, \lambda _k) \}_{k=0}^{N} \right)  = 0
	\end{equation}
	is equivalent to the discrete constrained right Hamilton's equations
	\begin{equation} \label{Discrete Constrained H equationsAppendix}
		q_{k+1}  = D_2 H_d^+(q_{k},p_{k+1}),  \qquad p_k  =  D_1H_d^+(q_k,p_{k+1}) + \langle \lambda_k, \nabla\mathcal{C}(q_k)  \rangle, \qquad  \mathcal{C}(q_k)  =0, 
	\end{equation}
	where $H_d^+ (q_k,p_{k+1}) $ is defined via equation \eqref{Discrete Right Hamiltonian}.
	\small \proof{Using the fact that $\delta q_0 =0$ and $\delta p_N =0$ since $(q_0,p_N)$ is fixed, we obtain the following expression for the variations of  $ \mathfrak{S}_d^{+} $:
	\begin{align*}
	\delta \mathfrak{S}_d^+  & = \delta \left(  p_N q_N - \sum_{k=0}^{N-1}{\left[ p_{k+1} q_{k+1} - H_d^+(q_k,p_{k+1}) - \langle \lambda_k, \mathcal{C}(q_k)  \rangle \right]} \right)  \\ &  = \delta \left(  -\sum_{k=0}^{N-2}{p_{k+1}q_{k+1}} + \sum_{k=0}^{N-1}{\left[  H_d^+(q_k,p_{k+1}) + \langle \lambda_k, \mathcal{C}(q_k)  \rangle \right]} \right)   \\
	& = -\sum_{k=0}^{N-2}{\left( q_{k+1}\delta p_{k+1} + p_{k+1}\delta q_{k+1} \right)}  + \sum_{k=0}^{N-1}{\left( D_1H_d^+(q_k,p_{k+1}) \delta q_k  + D_2 H_d^+(q_k,p_{k+1}) \delta p_{k+1} \right)}  \\ & \qquad  \qquad + \sum_{k=0}^{N-1}{\left( \langle \lambda_k ,\nabla \mathcal{C}(q_k) \delta q_k \rangle    + \langle \delta \lambda_k , \mathcal{C}(q_k) \rangle  \right)} \\ 
	& = -\sum_{k=1}^{N-1}{\left( q_{k}\delta p_{k} + p_{k}\delta q_{k} \right)} + \sum_{k=1}^{N-1}{ D_1H_d^+(q_k,p_{k+1}) \delta q_k }  +  \sum_{k=0}^{N-2}{D_2 H_d^+(q_k,p_{k+1}) \delta p_{k+1}} \\ & \qquad  \qquad + \sum_{k=0}^{N-1}{\left( \langle \lambda_k ,\nabla \mathcal{C}(q_k) \delta q_k \rangle    + \langle \delta \lambda_k , \mathcal{C}(q_k) \rangle  \right)} \\  
	& = -\sum_{k=1}^{N-1}{\left( q_{k}\delta p_{k} + p_{k}\delta q_{k} \right)} + \sum_{k=1}^{N-1}{ D_1H_d^+(q_k,p_{k+1}) \delta q_k }  +   \sum_{k=1}^{N-1}{D_2 H_d^+(q_{k-1},p_{k}) \delta p_{k}} \\ & \qquad \qquad  + \sum_{k=0}^{N-1}{\left( \langle \lambda_k ,\nabla \mathcal{C}(q_k) \delta q_k \rangle    + \langle \delta \lambda_k , \mathcal{C}(q_k) \rangle  \right)} \\ & =   \sum_{k=1}^{N-1}{  \left[  -q_k  +D_2 H_d^+(q_{k-1},p_{k})    \right] \delta p_k } +   \sum_{k=0}^{N-1}{ \langle \delta \lambda_k , \mathcal{C}(q_k) \rangle  }  +  \sum_{k=1}^{N-1}{  \left[  -p_k  + D_1H_d^+(q_k,p_{k+1}) + \langle \lambda_k, \nabla  \mathcal{C}(q_k)  \rangle \right] \delta q_k }. 	
\end{align*}
If the discrete constrained right Hamilton's equations \eqref{Discrete Constrained H equationsAppendix} are satisfied, then each term vanishes and $\delta \mathfrak{S}_d^{+}= 0$. Conversely, if $\delta \mathfrak{S}_d^{+} = 0$, then a discrete fundamental theorem of the calculus of variations yields the discrete constrained right Hamilton's equations \eqref{Discrete Constrained H equationsAppendix}. \qed \\
}
\end{theorem}

\subsection{Proof of Theorem~\ref{theorem: Discrete Variational Principle III}} \label{Appendix: Type III Discrete Variational}

	\begin{theorem} 
	The Type III discrete Hamilton's phase space variational principle
	\begin{equation}
		\delta \mathfrak{S}_d^{-} \left(\{  (q_k,p_k, \lambda _k) \}_{k=0}^{N} \right)  = 0
	\end{equation}
	is equivalent to the discrete constrained left Hamilton's equations
	\begin{equation} \label{Discrete Constrained H equationsAppendix III}
		q_{k}  = -D_2 H_d^-(q_{k+1},p_{k}),  \qquad p_{k+1}  =  -D_1H_d^-(q_{k+1},p_{k}) - \langle \lambda_{k+1}, \nabla\mathcal{C}(q_{k+1})  \rangle, \qquad  \mathcal{C}(q_k)  =0, 
	\end{equation}
	where $H_d^- (q_{k+1},p_{k}) $ is defined via equation \eqref{Discrete Left Hamiltonian}.
\small 	\proof{Using the fact that $\delta q_N =0$ and $\delta p_0 =0$ since $(q_N,p_0)$ is fixed, we obtain the following expression for the variations of  $ \mathfrak{S}_d^{-} $:
		\begin{align*}
			\delta \mathfrak{S}_d^-  & = \delta \left(  - p_0 q_0 - \sum_{k=0}^{N-1}{\left[ -p_{k} q_{k} - H_d^-(q_{k+1},p_{k}) - \langle \lambda_{k+1}, \mathcal{C}(q_{k+1})  \rangle  \right]} \right)  \\ &  = \delta \left(  \sum_{k=1}^{N-1}{p_{k}q_{k}} + \sum_{k=0}^{N-1}{\left[  H_d^-(q_{k+1},p_{k}) + \langle \lambda_{k+1}, \mathcal{C}(q_{k+1})  \rangle \right]} \right)   \\
			& = \sum_{k=1}^{N-1}{\left( q_{k}\delta p_{k} + p_{k}\delta q_{k} \right)}  + \sum_{k=0}^{N-1}{\left( D_1H_d^-(q_{k+1},p_{k}) \delta q_{k+1}  + D_2 H_d^-(q_{k+1},p_{k}) \delta p_{k} \right)}   \\ & \qquad \qquad   + \sum_{k=0}^{N-1}{\left( \langle \lambda_{k+1} ,\nabla \mathcal{C}(q_{k+1}) \delta q_{k+1} \rangle    + \langle \delta \lambda_{k+1} , \mathcal{C}(q_{k+1}) \rangle  \right)} \\ 
			& =  \sum_{k=0}^{N}{\left( q_{k}\delta p_{k} + p_{k}\delta q_{k} \right)} + \sum_{k=0}^{N-2}{ D_1H_d^-(q_{k+1},p_{k}) \delta q_{k+1} }  +  \sum_{k=1}^{N-1}{D_2 H_d^-(q_{k+1},p_{k}) \delta p_{k}}  \\ & \qquad \qquad  + \sum_{k=0}^{N-1}{\left( \langle \lambda_{k+1} ,\nabla \mathcal{C}(q_{k+1}) \delta q_{k+1} \rangle    + \langle \delta \lambda_{k+1} , \mathcal{C}(q_{k+1}) \rangle  \right)} \\ 
			& =  \sum_{k=1}^{N-1}{\left( q_{k}\delta p_{k} + p_{k}\delta q_{k} \right)} + \sum_{k=1}^{N-1}{ D_1H_d^-(q_k,p_{k-1}) \delta q_k }  +  \sum_{k=1}^{N-1}{D_2 H_d^-(q_{k+1},p_{k}) \delta p_{k}}  \\ & \qquad \qquad   + \sum_{k=1}^{N}{\langle \lambda_{k} ,\nabla \mathcal{C}(q_{k}) \delta q_{k} \rangle  } +   \sum_{k=0}^{N-1}{ \langle \delta \lambda_{k+1} , \mathcal{C}(q_{k+1}) \rangle  }\\ & =   \sum_{k=1}^{N-1}{  \left[  q_k  + D_2 H_d^-(q_{k+1},p_{k})    \right] \delta p_k } +  \sum_{k=0}^{N-1}{ \langle \delta \lambda_{k+1} , \mathcal{C}(q_{k+1}) \rangle}  +  \sum_{k=1}^{N-1}{  \left[   p_k  + D_1H_d^-(q_k,p_{k-1}) + \langle \lambda_k ,\nabla \mathcal{C}(q_k)  \rangle  \right] \delta q_k }.
		\end{align*}
		If the discrete constrained left Hamilton's equations \eqref{Discrete Constrained H equationsAppendix III} are satisfied, then each term vanishes and $\delta \mathfrak{S}_d^{-} = 0$. Conversely, if $\delta \mathfrak{S}_d^{-} = 0$, then a discrete fundamental theorem of the calculus of variations yields the discrete constrained left Hamilton's equations \eqref{Discrete Constrained H equationsAppendix III}. \qed \\
	}
\end{theorem}

\bibliography{ConstrainedBib}

\def\cprime{$'$}
\begin{thebibliography}{48}
\providecommand{\natexlab}[1]{#1}
\providecommand{\url}[1]{\texttt{#1}}
\expandafter\ifx\csname urlstyle\endcsname\relax
  \providecommand{\doi}[1]{doi: #1}\else
  \providecommand{\doi}{doi: \begingroup \urlstyle{rm}\Url}\fi

\bibitem[Abraham et~al.(1988)Abraham, Marsden, and Ratiu]{AbMaRa1988}
R.~Abraham, J.~E. Marsden, and T.~Ratiu.
\newblock \emph{Manifolds, Tensor Analysis, and Applications}, volume~75 of
  \emph{Applied Mathematical Sciences}.
\newblock Springer, New York, second edition, 1988.

\bibitem[Absil et~al.(2008)Absil, Mahony, and Sepulchre]{Absil2008}
P.~A. Absil, R.~Mahony, and R.~Sepulchre.
\newblock \emph{Optimization Algorithms on Matrix Manifolds}.
\newblock Princeton University Press, Princeton, NJ, 2008.
\newblock ISBN 978-0-691-13298-3.

\bibitem[Ahn and Sra(2020)]{Sra2020}
K.~Ahn and S.~Sra.
\newblock From {N}esterov's estimate sequence to {R}iemannian acceleration.
\newblock In \emph{Proceedings of Thirty Third Conference on Learning Theory},
  volume 125 of \emph{Proceedings of Machine Learning Research}, pages 84--118.
  PMLR, 09--12 Jul 2020.

\bibitem[Alimisis et~al.(2020{\natexlab{a}})Alimisis, Orvieto, B\'ecigneul, and
  Lucchi]{Alimisis2020-1}
F.~Alimisis, A.~Orvieto, G.~B\'ecigneul, and A.~Lucchi.
\newblock Practical accelerated optimization on {R}iemannian manifolds,
  2020{\natexlab{a}}.

\bibitem[Alimisis et~al.(2020{\natexlab{b}})Alimisis, Orvieto, B\'ecigneul, and
  Lucchi]{alimisis2020}
F.~Alimisis, A.~Orvieto, G.~B\'ecigneul, and A.~Lucchi.
\newblock A continuous-time perspective for modeling acceleration in
  {R}iemannian optimization.
\newblock In \emph{Proceedings of the 23rd International AISTATS Conference},
  volume 108 of \emph{PMLR}, pages 1297--1307, 2020{\natexlab{b}}.

\bibitem[Alimisis et~al.(2021)Alimisis, Orvieto, B{\'e}cigneul, and
  Lucchi]{Alimisis2021}
F.~Alimisis, A.~Orvieto, G.~B{\'e}cigneul, and A.~Lucchi.
\newblock Momentum improves optimization on {R}iemannian manifolds.
\newblock In \emph{AISTATS}, 2021.

\bibitem[Arnol{\cprime}d(1989)]{Ar1989}
V.~I. Arnol{\cprime}d.
\newblock \emph{Mathematical methods of classical mechanics}, volume~60 of
  \emph{Graduate Texts in Mathematics}.
\newblock Springer-Verlag, New York, second edition, 1989.
\newblock Translated from the Russian by K. Vogtmann and A. Weinstein.

\bibitem[Benettin and Giorgilli(1994)]{Benettin1994}
G.~Benettin and A.~Giorgilli.
\newblock On the {H}amiltonian interpolation of near-to-the identity symplectic
  mappings with application to symplectic integration algorithms.
\newblock \emph{Journal of Statistical Physics}, 74:\penalty0 1117--1143, 03
  1994.
\newblock \doi{10.1007/BF02188219}.

\bibitem[Benito and Mart\'in~de Diego(2005)]{Benito2005}
R.~Benito and D.~Mart\'in~de Diego.
\newblock Discrete vakonomic mechanics.
\newblock \emph{Journal of Mathematical Physics}, 46\penalty0 (8):\penalty0
  083521, 2005.
\newblock \doi{10.1063/1.2008214}.

\bibitem[Boumal(2020)]{Boumal2020}
N.~Boumal.
\newblock An introduction to optimization on smooth manifolds, 2020.
\newblock URL \url{http://www.nicolasboumal.net/book}.

\bibitem[Cort{\'e}s and Mart{\'i}nez(2001)]{Cortes2001}
J.~Cort{\'e}s and S.~Mart{\'i}nez.
\newblock Non-holonomic integrators.
\newblock \emph{Nonlinearity}, 14:\penalty0 1365--1392, 2001.

\bibitem[Cort\'es et~al.(2002)Cort\'es, de~Le\'on, Mart\'in~de Diego, and
  Martínez]{Cortes2002}
J.~Cort\'es, M.~de~Le\'on, D.~Mart\'in~de Diego, and S.~Martínez.
\newblock Geometric description of vakonomic and nonholonomic dynamics.
  comparison of solutions.
\newblock \emph{SIAM Journal on Control and Optimization}, 41\penalty0
  (5):\penalty0 1389--1412, 2002.
\newblock \doi{10.1137/S036301290036817X}.

\bibitem[{de Le\'on} et~al.(2004){de Le\'on}, Mart\'in~de Diego, and
  Santamaría-Merino]{Leon2003}
M.~{de Le\'on}, D.~Mart\'in~de Diego, and A.~Santamaría-Merino.
\newblock Geometric numerical integration of nonholonomic systems and optimal
  control problems.
\newblock \emph{European Journal of Control}, 10\penalty0 (5):\penalty0
  515--521, 2004.

\bibitem[{de Le\'on} et~al.(2007){de Le\'on}, Mart\'in~de Diego, and
  Santamaría-Merino]{Leon2007}
M.~{de Le\'on}, D.~Mart\'in~de Diego, and A.~Santamaría-Merino.
\newblock Discrete variational integrators and optimal control theory.
\newblock \emph{Advances in Computational Mathematics}, 26:\penalty0 251--268,
  2007.
\newblock \doi{10.1007/s10444-004-4093-5}.

\bibitem[Duruisseaux and Leok(2021)]{duruisseaux2021projection}
V.~Duruisseaux and M.~Leok.
\newblock Accelerated optimization on {R}iemannian manifolds via projected
  variational integrators.
\newblock 2021.

\bibitem[Duruisseaux and Leok(2022{\natexlab{a}})]{Duruisseaux2021Riemannian}
V.~Duruisseaux and M.~Leok.
\newblock A variational formulation of accelerated optimization on {R}iemannian
  manifolds.
\newblock \emph{SIAM Journal on Mathematics of Data Science},
  2022{\natexlab{a}}.
\newblock accepted.

\bibitem[Duruisseaux and Leok(2022{\natexlab{b}})]{Duruisseaux2022Lagrangian}
V.~Duruisseaux and M.~Leok.
\newblock Time-adaptive {L}agrangian variational inegrators for accelerated
  optimization on manifolds,.
\newblock 2022{\natexlab{b}}.
\newblock URL \url{https://arxiv.org/abs/2201.03774}.

\bibitem[Duruisseaux et~al.(2021)Duruisseaux, Schmitt, and
  Leok]{duruisseaux2020adaptive}
V.~Duruisseaux, J.~Schmitt, and M.~Leok.
\newblock Adaptive {H}amiltonian variational integrators and applications to
  symplectic accelerated optimization.
\newblock \emph{SIAM Journal on Scientific Computing}, 43\penalty0
  (4):\penalty0 A2949--A2980, 2021.
\newblock URL \url{https://doi.org/10.1137/20M1383835}.

\bibitem[Eld{\'{e}}n and Park(1999)]{Elden1999}
L.~Eld{\'{e}}n and H.~Park.
\newblock A {P}rocrustes problem on the {S}tiefel manifold.
\newblock \emph{Numerische Mathematik}, 82\penalty0 (4):\penalty0 599--619,
  1999.
\newblock \doi{10.1007/s002110050432}.

\bibitem[Fedorov and Zenkov(2005)]{Fedorov2005}
Y.~N. Fedorov and D.~V. Zenkov.
\newblock Discrete nonholonomic {LL} systems on {L}ie groups.
\newblock \emph{Nonlinearity}, 18:\penalty0 2211--2241, 2005.

\bibitem[Golub and Van~Loan(2013)]{Golub2013}
G.~H. Golub and C.~F. Van~Loan.
\newblock \emph{Matrix Computations}.
\newblock Johns Hopkins Studies in the Mathematical Sciences. Johns Hopkins
  University Press, 2013.
\newblock ISBN 9781421407944.

\bibitem[Hairer et~al.(2006)Hairer, Lubich, and Wanner]{HaLuWa2006}
E.~Hairer, C.~Lubich, and G.~Wanner.
\newblock \emph{Geometric {N}umerical {I}ntegration}, volume~31 of
  \emph{Springer Series in Computational Mathematics}.
\newblock Springer-Verlag, Berlin, second edition, 2006.

\bibitem[{Hall} and {Leok}(2015)]{HaLe2012}
J.~{Hall} and M.~{Leok}.
\newblock {Spectral Variational Integrators}.
\newblock \emph{Numer. Math.}, 130\penalty0 (4):\penalty0 681--740, 2015.

\bibitem[Holm et~al.(2009)Holm, Schmah, and Stoica]{Holm2009}
D.~Holm, T.~Schmah, and C.~Stoica.
\newblock \emph{Geometric Mechanics and Symmetry: From Finite to Infinite
  Dimensions}.
\newblock Oxford Texts in Applied and Engineering Mathematics. OUP Oxford,
  2009.
\newblock ISBN 9780199212910.

\bibitem[Jim{\'e}nez and Mart\'in~de Diego(2012)]{Jimenez2012}
F.~Jim{\'e}nez and D.~Mart\'in~de Diego.
\newblock Continuous and discrete approaches to vakonomic mechanics.
\newblock \emph{Revista De La Real Academia De Ciencias Exactas Fisicas Y
  Naturales Serie A-matematicas}, 106, 03 2012.
\newblock \doi{10.1007/s13398-011-0028-4}.

\bibitem[Jost(2017)]{Jost2017}
J.~Jost.
\newblock \emph{Riemannian {G}eometry and {G}eometric {A}nalysis}.
\newblock Universitext. Springer, Cham, 7th edition, 2017.

\bibitem[Junge et~al.(2005)Junge, Marsden, and Ober-Bl\"obaum]{Junge2005}
O.~Junge, J.~E. Marsden, and S.~Ober-Bl\"obaum.
\newblock Discrete mechanics and optimal control.
\newblock \emph{IFAC Proceedings Volumes}, 38\penalty0 (1):\penalty0 538--543,
  2005.

\bibitem[Lang(1999)]{Lang1999}
S.~Lang.
\newblock \emph{Fundamentals of Differential Geometry}, volume 191 of
  \emph{Graduate Texts in Mathematics}.
\newblock Springer -Verlag, New York, 1999.
\newblock ISBN 9780387985930.

\bibitem[Lee(2018)]{Lee2019}
{J. M.} Lee.
\newblock \emph{Introduction to {R}iemannian Manifolds}, volume 176 of
  \emph{Graduate Texts in Mathematics}.
\newblock Springer, Cham, second edition, 2018.

\bibitem[Leok and Ohsawa(2011)]{LeOh2008}
M.~Leok and T.~Ohsawa.
\newblock Variational and geometric structures of discrete {D}irac mechanics.
\newblock \emph{Found. Comput. Math.}, 11\penalty0 (5):\penalty0 529--562,
  2011.

\bibitem[Leok and Zhang(2011)]{LeZh2011}
M.~Leok and J.~Zhang.
\newblock Discrete {H}amiltonian variational integrators.
\newblock \emph{IMA Journal of Numerical Analysis}, 31\penalty0 (4):\penalty0
  1497--1532, 2011.

\bibitem[Liu et~al.(2017)Liu, Shang, Cheng, Cheng, and Jiao]{Liu2017}
Y.~Liu, F.~Shang, J.~Cheng, H.~Cheng, and L.~Jiao.
\newblock Accelerated first-order methods for geodesically convex optimization
  on {R}iemannian manifolds.
\newblock In \emph{NeurIPS}, volume~30, pages 4868--4877, 2017.

\bibitem[Marsden and Ratiu(1999)]{MaRa1999}
{J. E.} Marsden and {T. S.} Ratiu.
\newblock \emph{Introduction to mechanics and symmetry}, volume~17 of
  \emph{Texts in Applied Mathematics}.
\newblock Springer-Verlag, New York, second edition, 1999.

\bibitem[Marsden and West(2001)]{MaWe2001}
J.~E. Marsden and M.~West.
\newblock Discrete mechanics and variational integrators.
\newblock \emph{Acta Numer.}, 10:\penalty0 357--514, 2001.

\bibitem[McLachlan and Perlmutter(2006)]{McLachlan2006}
R.~I. McLachlan and M.~Perlmutter.
\newblock Integrators for nonholonomic mechanical systems.
\newblock \emph{Journal of Nonlinear Science}, 16:\penalty0 283--328, 2006.

\bibitem[Nash(1956)]{Nash1956}
J.~Nash.
\newblock The imbedding problem for {R}iemannian manifolds.
\newblock \emph{Annals of Mathematics}, 63\penalty0 (1):\penalty0 20--63, 1956.
\newblock ISSN 0003486X.

\bibitem[Nesterov(1983)]{Nes83}
Y.~Nesterov.
\newblock A method of solving a convex programming problem with convergence
  rate $\mathcal{O}(1/k^2)$.
\newblock \emph{Soviet Mathematics Doklady}, 27\penalty0 (2):\penalty0
  372--376, 1983.

\bibitem[Nesterov(2004)]{Nes04}
Y.~Nesterov.
\newblock \emph{Introductory Lectures on Convex Optimization: A Basic Course},
  volume~87 of \emph{Applied Optimization}.
\newblock Kluwer Academic Publishers, Boston, MA, 2004.

\bibitem[Reich(1999)]{Re1999}
S.~Reich.
\newblock Backward error analysis for numerical integrators.
\newblock \emph{SIAM J. Numer. Anal.}, 36:\penalty0 1549--1570, 1999.

\bibitem[Schmitt and Leok(2017)]{ScLe2017}
J.~M. Schmitt and M.~Leok.
\newblock Properties of {H}amiltonian variational integrators.
\newblock \emph{IMA Journal of Numerical Analysis}, 38\penalty0 (1):\penalty0
  377--398, 03 2017.

\bibitem[Schmitt et~al.(2018)Schmitt, Shingel, and Leok]{ScShLe2017}
J.~M. Schmitt, T.~Shingel, and M.~Leok.
\newblock {L}agrangian and {H}amiltonian {T}aylor variational integrators.
\newblock \emph{{BIT} Numerical Mathematics}, 58:\penalty0 457--488, 2018.
\newblock \doi{10.1007/s10543-017-0690-9}.

\bibitem[Su et~al.(2016)Su, Boyd, and Candes]{SuBoCa16}
W.~Su, S.~Boyd, and E.~Candes.
\newblock A differential equation for modeling {N}esterov's {A}ccelerated
  {G}radient method: theory and insights.
\newblock \emph{Journal of Machine Learning Research}, 17\penalty0
  (153):\penalty0 1--43, 2016.

\bibitem[Sutskever et~al.(2013)Sutskever, Martens, Dahl, and
  Hinton]{Sutskever2013}
I.~Sutskever, J.~Martens, G.~Dahl, and G.~Hinton.
\newblock On the importance of initialization and momentum in deep learning.
\newblock In \emph{Proceedings of the 30th International Conference on
  International Conference on Machine Learning - Volume 28}, ICML'13, pages
  1139--1147, Atlanta, GA, USA, 2013.

\bibitem[Whitney(1944{\natexlab{a}})]{Whitney1944_1}
H.~Whitney.
\newblock The singularities of a smooth $n$-manifold in $(2n - 1)$-space.
\newblock \emph{Annals of Mathematics}, 45\penalty0 (2):\penalty0 247--293,
  1944{\natexlab{a}}.
\newblock ISSN 0003486X.

\bibitem[Whitney(1944{\natexlab{b}})]{Whitney1944_2}
H.~Whitney.
\newblock The self-intersections of a smooth $n$-manifold in $2n$-space.
\newblock \emph{Annals of Mathematics}, 45\penalty0 (2):\penalty0 220--246,
  1944{\natexlab{b}}.
\newblock ISSN 0003486X.

\bibitem[Wibisono et~al.(2016)Wibisono, Wilson, and Jordan]{WiWiJo16}
A.~Wibisono, A.~Wilson, and M.~Jordan.
\newblock A variational perspective on accelerated methods in optimization.
\newblock \emph{Proceedings of the National Academy of Sciences}, 113\penalty0
  (47):\penalty0 E7351--E7358, 2016.

\bibitem[Zhang and Sra(2016)]{Sra2016}
H.~Zhang and S.~Sra.
\newblock First-order methods for geodesically convex optimization.
\newblock In \emph{29th Annual Conference on Learning Theory}, pages
  1617--1638, 2016.

\bibitem[Zhang and Sra(2018)]{Sra2018}
H.~Zhang and S.~Sra.
\newblock An estimate sequence for geodesically convex optimization.
\newblock In \emph{Proceedings of the 31st Conference On Learning Theory},
  volume~75 of \emph{Proceedings of Machine Learning Research}, pages
  1703--1723, Jul 2018.

\end{thebibliography}
\bibliographystyle{plainnat}

\end{document}